\documentclass[11pt,fleqn]{article}
\usepackage{latexsym}
\usepackage{graphics}
\usepackage{graphicx}
\usepackage{amsmath}
\usepackage{amsfonts}
\usepackage{eucal}
\usepackage{setspace}
\usepackage{titlesec}

\usepackage[left=1in,top=0.9in,right=1in,bottom=1.4in,nohead]{geometry}

\newcommand{\ind}{\setlength{\parindent}{2em}}

\def\h{\hspace{.5em}}

\def\v{\vspace{.5ex}}
\def\vv{\vspace{1ex}}
\def\vvv{\vspace{2ex}}

\def\h{\hspace{.5em}}

\titleformat{\section} 
{\normalfont\bfseries}
{\thesection.}{0.5em}{}

\title{Determining isotopy classes of crossing arcs in alternating links}
\author{Anastasiia Tsvietkova}
\date{}

\begin{document}
\everymath{\displaystyle}
\maketitle
 \footnotesize
\begin{description} \item \qquad \textbf{ Abstract.} Given a reduced alternating diagram for a link, we obtain conditions that guarantee that the link complement has a complete hyperbolic structure, crossing arcs are the edges of an ideal geodesic triangulation, and every crossing arc is isotopic to a simple geodesic. The latter was conjectured by Sakuma and Weeks in 1995. We provide infinite families of closed braids for which our conditions hold. 

\vvv


\vvv


\end{description}

\normalsize

\section{Overview}

A link diagram provides a combinatorial description of a topological object, a link complement in $S^3$. A natural
question arising from this description is whether the complement
can be endowed with a complete hyperbolic structure. This question is connected with the question of the existence of an ideal geodesic triangulation of the complement. While both of these questions can often be answered after tedious computations for a particular link, it is not \textit{a priori} clear what would be a successful choice of edges for such a triangulation.

In this note, we show that under certain conditions (that can easily be checked) the crossing arcs of a reduced alternating diagram are the edges of an ideal geodesic triangulation. By a crossing arc we mean a cusp-to-cusp arc traveling from the underpass to the overpass of a crossing. The triangulations we obtain induce the complete hyperbolic metric on the link complement (which implies, in particular, that the link is hyperbolic without a reference to Geometrization). It follows that, under these conditions, crossing arcs are isotopic to simple geodesics. We provide examples for which this holds, including an infinite family of links.

This question has a long history. W. Thurston noticed that if we choose a decomposition of a hyperbolic link complement into two polyhedra, where the edges are crossing arcs of a link diagram, then we often obtain an ideal geodesic triangulation from it by subdividing the polyhedra (\cite{Thurston}). The method was generalized by Menasco for alternating links (\cite{Menasco}), and by Petronio beyond alternating (\cite{Carlo}). While this suggests that the arcs are often isotopic to geodesics, the procedure fails to provide an ideal geodesic triangulation with positive-volume tetrahedra in general. On the other hand, Sakuma and Weeks conjectured that crossing arcs of a reduced alternating diagram are the arcs of the canonical cell decomposition of the link complement (\cite{SakumaWeeks}), which would imply they are isotopic to geodesics. Their conjecture was proved for hyperbolic 2-bridge links in \cite{Akiyoshi} and independently in \cite{Gueritad} (see also Appendix to \cite{FuterAppendix}), but has not been established for any wider classes of links. 

Beyond 2-bridge links, there were few additional results that identify cusp-to-cusp arcs in link complements with certain topological or combinatorial description as geodesics. Sakuma and Weeks gave examples of canonical cell decompositions implying that crossing arcs of some families of symmetric alternating links are isotopic to geodesics (see Examples I.2.2-I.2.4 in \cite{SakumaWeeks}). Work of Adams, Burton, Cooper, Futer, Purcell provided information about isotopy classes of certain tunnel arcs under additional restrictions (\cite{Adams, AdamsReid, Burton, Cooper}).

The sufficient conditions we obtain turn out to have a simple geometric interpretation. Let us return to the polyhedral decomposition suggested by Thurston and described by Menasco for alternating links. The term ``polyhedron" is used only in topological sense here, since even for a hyperbolic link the faces might not be planar, \textit{i.e.} might not lie in one hyperbolic plane. We will call such a polyhedron \textit{cross-sectionally convex}, if at every ideal vertex all interior angles of a cross-section are in $(0,\pi)$. Note that this is not equivalent to the usual convexity: the faces might still be non-planar and, when subdivided, might yield a non-convex polyhedron.

 Our conditions imply cross-sectional convexity of the two polyhedra suggested by Thurston. Then we show that every cross-sectionally convex polyhedral decomposition can be subdivided into an ideal geodesic \textit{partially flat} triangulation (``partially flat" means that some, but not all tetrahedra have 0 volume, while the rest have positive volume). In particular, for links satisfying our conditions the purely combinatorial algorithm described by Thurston, Menasco and Petronio provides an ideal geodesic partially flat triangulation that induces the complete hyperbolic structure on the link complement. Such triangulations appear to be useful for various other purposes as well (see, for example, \cite{LackenbyHeegaard} or \S 6 in \cite{NeumannZagier}).

The following is a short overview of our methods and techniques. In \cite{method, thesis}, Morwen Thistlethwaite and the author introduced an alternative way of parameterizing hyperbolic structures of links.  It uses complex labels assigned to a link diagram that describe horoball structures in $\mathbb{H}^3$. We will refer to these labels as to \textit{diagram labels}. The triangulation is not performed, and instead the method uses isometries of preimages of polygons bounded by the regions of the link diagram. This results in a set of relations, to which we will refer as to \textit{hyperbolicity relations}. 

In \cite{method}, we start with a hyperbolic link, and then describe the relations for the diagram labels merely as a convenient method for computing the (already existent) complete hyperbolic structure of the complement. In this paper instead we start with an arbitrary link diagram and the complex labels that satisfy the hyperbolicity relations for this diagram, and then establish additional conditions on the labels that guarantee the existence of the induced complete hyperbolic structure.

 The paper is organized as follows. In Section 2, we describe the setting, developing a model that gives the complete hyperbolic structure based on diagram labels rather than a triangulation process. In Section 3, we lay out and explain the conditions on diagram labels sufficient for our purposes, and prove that under these conditions, the Thurston-Menasco polyhedral decomposition, with faces triangulated, is properly embedded in $\mathbb{H}^3$. In Section 4, we show that, further, the decomposition yields an ideal partially flat geodesic triangulation, and use this to conclude that a link complement has a complete hyperbolic structure. In section 5, we prove that crossing arcs are isotopic to simple geodesics under our conditions. Section 6 gives examples for which this holds (including infinite families of links), demonstrating that one can check the conditions from a link diagram. The families of links are different from the links provided by Sakuma and Weeks in \cite{SakumaWeeks}.

\section{Diagram labels}

Consider the complement of a link $L$ in $S^3$, and a diagram $D$ of $L$. In what follows, we will describe the correspondence between the points on the peripheral boundary of $S^3-L$ and the points in the boundary of $\mathbb{H}^3$ through a solution of the hyperbolicity equations, introduced in \cite{method}. In further sections, we will investigate under which conditions on the labels this correspondence yields the complete hyperbolic structure and the developing map for the manifold $S^3-L$.

 Every region $R$ of $D$ that is incident to at least three crossings may be viewed as a disk bounded by the geodesic arcs traveling from the overpass to the underpass at every crossing of $R$, and by the arcs traveling on the boundary torus from one crossing of $R$ to the next crossing of $R$ (black and gray arcs on Fig.1, left, respectively, where a region incident to three crossings is depicted). Every arc can be assigned a complex number which we call a diagram label. When the complete hyperbolic structure exists on $S^3-L$, the diagram labels are called crossing and edge labels respectively, and they parameterize the hyperbolic structure (see \cite{method} for details and geometric definitions of the labels). 
 
 \begin{center}
\includegraphics[scale=0.32]{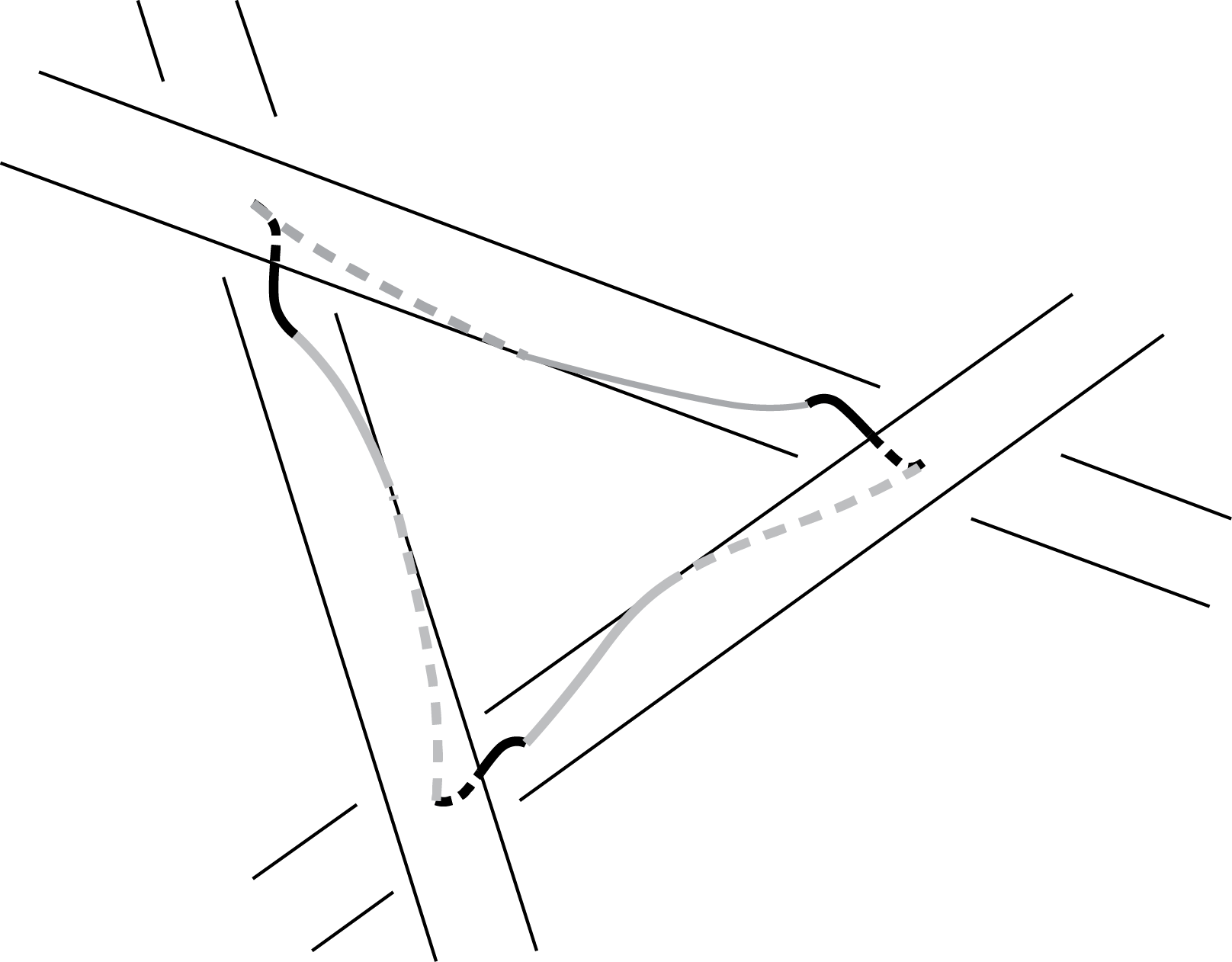} \hspace{0.5in}
\includegraphics[scale=0.38]{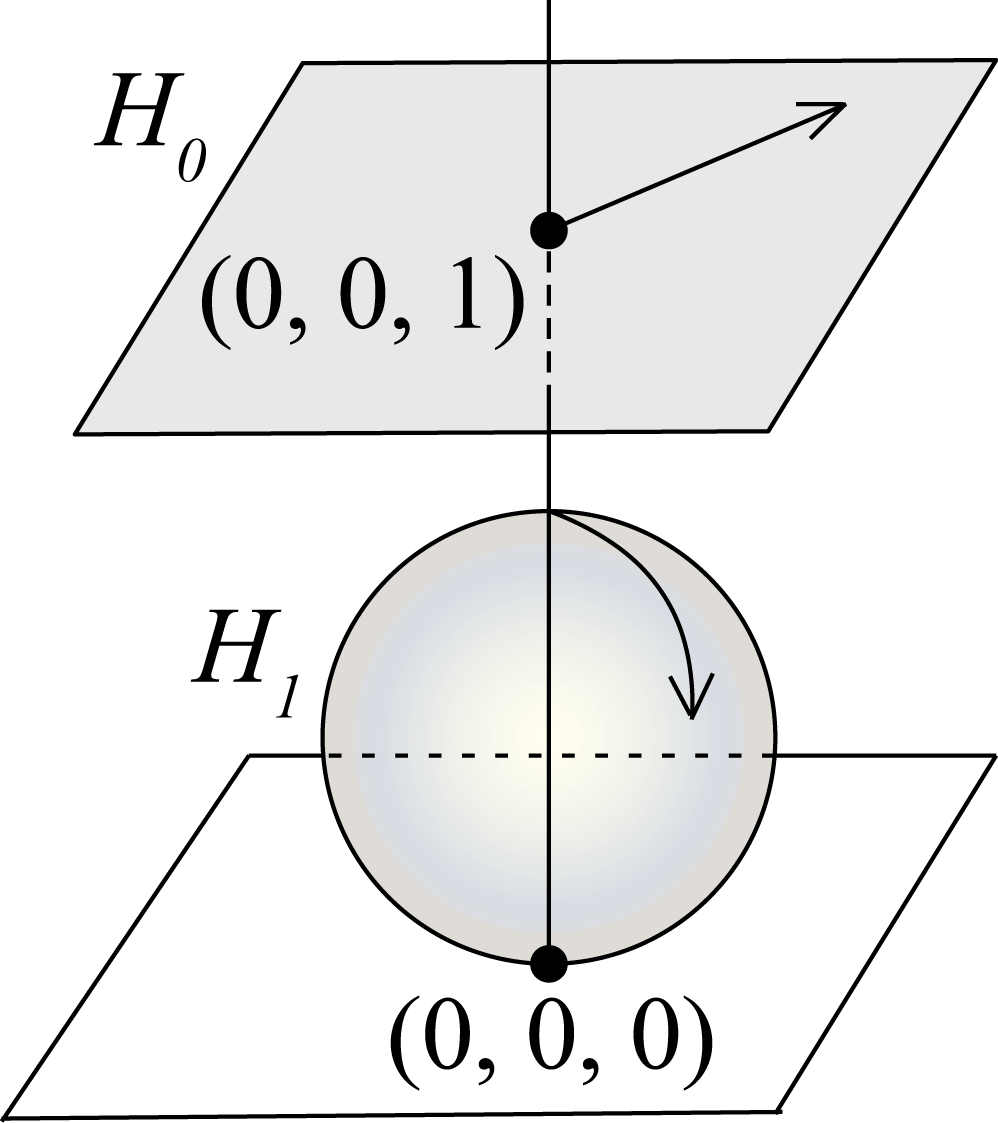} 

 Fig.1
 \end{center}

 Hyperbolicity relations are two sets of polynomial relations in the diagram labels. The first set, called the region relations, guarantees that the composition of hyperbolic isometries rotating a preimage in $\mathbb{H}^3$ of the boundary of each disk, corresponding to a region of the diagram, must be the identity. If the region is 2-sided (\textit{i.e.} it is a bigon), the two edge labels inside it are set to 0, and the two crossing labels are set equal, which corresponds to having no disk, but rather two homotopic geodesic arcs at crossings. The second set of relations (called the edge relations) for alternating links consists of relations of the form either $u=v\pm1$ or $u=v$, where $u,v$ are two edge labels assigned to two different sides of the same edge of $D$. By an edge of $D$ we mean a segment from a crossing to the next crossing, and hence $u$ and $v$ lie in two different regions of $D$ separated by this segment.  The edge relations guarantee that whenever two arcs on the boundary torus form a simple closed curve going around a strand of $L$, the curve is homotopic to a meridian of length 1. We refer the reader to \cite{method} for further details and examples.
 
 Assume that there is a complex solution $\overline{x}$ to the hyperbolicity relations  for \( D \) (it is possible that all entries of $\overline{x}$ have 0 imaginary part, and therefore are real). In particular, if there are $n$ crossings in $D$,  $\overline{x}$ consists of $n$ crossing labels $w_1, w_2, ..., w_n$, and $2n$ edge labels $u_1, u_2, ..., u_{2n}$.

Locally, every cusp (a neighborhood of a boundary torus for every link component) can be endowed with a Euclidean structure.  At an overpass or an underpass of a crossing of $D$, choose a cusp cross-section with a unit meridian. We will now describe the correspondence between the link complement and a picture that, for a hyperbolic link complement, will later prove to be its preimage in $\mathbb{H}^3$. Let the cusp cross-section correspond to an infinite union of Euclidean planes in this picture, and let the unit meridian correspond to the real number 1 on each of them. We will view each of these planes in the Euclidean three-dimensional space as a sphere minus the South pole, touching the plane $z=0$ from above. We may consider these spheres as horospheres in $\mathbb{H}^3$, using the upper half-space model of $\mathbb{H}^3=\{(x,y,z) \hspace{0.05cm} | \hspace{0.05 cm} z>0\}$.

In what follows, we will connect the points where the horospheres are tangent to the plane $z=0$ (we call such a point $P_i$ the center of the corresponding horosphere $H_i$) by hyperbolic geodesics. We will situate and scale the horospheres so that the geometric definitions of the diagram labels are satisfied, even though the link $L$ might not be hyperbolic. For this, we use a correspondence between the points $P_i$ on the boundary of $\mathbb{H}^3$ and the points $\bar{P}_i$ on the boundary torus of $S^3-L$. In particular, for each overpass or underpass of $D$ we consider one point $\bar{P}_i$ on the boundary torus, and a horosphere centered at a corresponding point $P_i$ on the boundary of $\mathbb{H}^3$.

To specify the size and the exact location of the horospheres, consider a region $R_0$ of $D$ that has at least three crossings (and hence $R_0$ is not a bigon). Consider
two consecutive crossings. The overstrand of the first crossing corresponds to a
point $P_0$ on the boundary of $\mathbb{H}^3$, and a horosphere $H_0$ centered at $P_0$. A crossing arc runs to
an understrand, corresponding to the point $P_1$ on the boundary of $\mathbb{H}^3$, and horosphere $H_1$. The
overstrand of the next crossing corresponds to the same point $P_1$ on $\mathbb{H}^3$, and the
same horosphere. A crossing arc runs from that overstrand to an understrand corresponding
to $P_2$ on the boundary of $\mathbb{H}^3$, and horosphere $H_2$. In $\mathbb{H}^3$, connect points $P_0$ and $P_1$, centers of $H_0$ and $H_1$, by the hyperbolic geodesic $\gamma_1$. Connect points $P_1$ and $P_2$ by a geodesic $\gamma_2$. Denote the preimage of the meridian on the horosphere $H_i$ by $m_i$ throughout. 

Let the crossing label that corresponds to the same crossing of $D$ that $\gamma_1$ does, be $w_1$ of $\bar{x}$. If $w_1$ is 0, then the centers of the horospheres $H_0, H_1$ coincide, $\gamma_1$ is null-homotopic, and we can proceed to the next horosphere. So let us assume that this crossing label is not 0.

Place one of the horospheres, $H_0$, as an ``infinite" horosphere, the plane $z=k$. Introduce the coordinates for $\mathbb{H}^3$ (\textit{i.e.} for the Euclidean 3-dimensional upper half-space $\{(x,y,z) \hspace{0.05cm} | \hspace{0.05 cm} z>0\}$) so that $k=1$, and the preimage of $m_0$ is the unit vector going from (0, 0, 1) to (1, 0, 1). Place $H_1$ at the point (0, 0, 0) (as on Fig.1, right, the horospheres are shaded). We may assume that $P_0, P_1$ are distinct, since otherwise $|w_1|=0$ (which happens if and only if $w_1=0$). Scale $H_1$ in $\mathbb{H}^3$ so that $|w_1|$ is the Euclidean diameter of $H_1$. Additionally, rotate $H_1$ so that the angle by rotation between $m_0$ and $m_1$ is $\arg{w_1}-\pi$. The horospheres $H_0, H_1$ are now situated so that the geometric definition of the crossing label $w_1$ from \cite{method} is satisfied.

The geodesic $\gamma_2$ is uniquely defined by one of its endpoints, $P_1$, and the corresponding edge label $u_2$ which tells where $\gamma_2$ pierces $H_1$. In particular, $m_1$, which is the preimage of a meridian on $H_1$, determines the unit distance and direction of the unit translation on $H_1$. The complex number $u_2$ (possibly, with a negative sign, depending on whether the direction of our travel along the region of $D$ corresponds to an orientation we choose for the link) is the translation on $H_1$ from the point where $\gamma_1$ pierces $H_1$ to the point (say, $M_1$) where $\gamma_2$ pierces $H_1$. Hence we obtain the position of $M_1$ on $H_1$. Two points, $P_1$ and $M_1$, uniquely determine the geodesic $\gamma_2$ and therefore determine the position of the other endpoint $P_2$ of $\gamma_2$.

The diameter of $H_2$ centered at $P_2$ and the direction of $m_2$ is determined by the next crossing label, $w_2$. Proceeding similarly from a horosphere to horosphere, region by region, we obtain a uniquely defined and scaled collection of horoballs with geodesics connecting them, a coordinate system on $\mathbb{H}^3$ and Euclidean coordinates on each horosphere such that the preimage of every meridian corresponds to the real number 1. 

Since the link $L$ is not necessarily a hyperbolic link, the described process might lead to certain degeneracy. For example, the picture may consist of a single horosphere, if all the labels are 0. In the next sections we will prove that this set-up together with a few additional conditions on the labels induce a complete hyperbolic metric on $S^3 - L$.

\section{Decomposition into two properly embedded polyhedra}

In this section, we consider a decomposition of the alternating link complement $S^3-L$ into two ideal polyhedra, $\overline{\Pi}_1, \overline{\Pi}_2$, one above the reduced alternating diagram $D$ of $L$, and one below, as described by Menasco in \cite{Menasco}. Assume additionally that the diagram $D$ is twist reduced in the sense of \cite{Lackenby}. Every alternating link admits such a diagram (see Section 3 in \cite{Lackenby} for the definition and explanation). 

The polyhedra are topological, and do not necessarily agree with the complete hyperbolic structure (which, possibly, does not even exist). However, the decomposition corresponds to the ideal complexes $\Pi_1, \Pi_2$ with straightened edges in $\mathbb{H}^3$ in the following way. In Section 2, we described the correspondence between points on the boundary of $\mathbb{H}^3$ and points on the boundary tori at overpasses/underpasses of $L$ using the diagram labels. This gives the correspondence between the vertices of $\bar{\Pi}_i$ and the vertices of $\Pi_i$, $i=1,2$. By ``straightened edges" we mean that every arc $\alpha$ of the decomposition of $S^3-L$ connecting two points on the boundary torus of $L$ (say, points $\bar{P}_1$ and $\bar{P}_2$) corresponds to a certain geodesic $\gamma$. The geodesic $\gamma$ connects the corresponding points $P_1, P_2$ in the boundary of $\mathbb{H}^3$.  Denote the set of such geodesic edges of $\Pi_i$ by $\mathcal{E}_i$.

Remove the vertex of $\Pi_i$ situated at infinity, and all edges of $\mathcal{E}_i$ incident to it, and denote the remaining collection of vertices and edges by $\mathcal{E}_i -\{v\}$. In the upper half-space model of $\mathbb{H}^3$,  let $f$ be the vertical projection of  $\mathcal{E}_i -\{v\}$ onto the plane $z=0$. A fragment of such a projection can be seen on Fig.4, right, with edges of the polyhedron in gray and their image in black.

Consider a collection of edges and vertices in $\mathbb{H}^3$ corresponding to a face $\bar{F}$ of $\bar{\Pi}_1$. Denote this collection by $\dot{F}$ and assume that $\dot{F}$ does not contain a vertex at infinity. The ideal vertices of $\dot{F}$ (denote all of them by $v_1, v_2, ..., v_j$) do not necessarily lie in one hyperbolic plane. We may however introduce a ``straightened'' face $F$ as follows. 

By taking the projection $f(\dot{F})$, we obtain  a polygon in the plane $z=0$. It is not necessarily a simple polygon, \textit{i.e.} its edges might intersect. There is however at least one simple polygon whose boundary includes at least one edge of $f(\dot{F})$. For every such simple polygon (say, $f(W_i')$), there is an ideal polygon $W'_i$ consisting of edges and ideal vertices of $\Pi_i$. Let the vertices of $W'_i$ be $v_1, ..., v_b$, and let the horospheres centered at them be $H_1, ...., H_b$ respectively. The size of the horospheres for a particular link is determined by the diagram labels, as described in the previous section. Now denote by $W_i$ the truncated polygon bounded by two types of arcs. An arc of the first type is a segment of a hyperbolic geodesic arc between a point on a horosphere $H_k$ (denote the point by $M_k$) and a point on a horosphere $H_{k+1}$ (denote the point by $N_{k+1}$), $k=1, 2,..., b-1$. Every such hyperbolic geodesic arc coincides with an edge of the ideal polygon $W_i'$. An arc of the second type is an arc traveling on a horosphere $H_k$ and connecting the points $N_k, M_k$, $k=1, 2, ..., b$. For every such truncated polygon $W_i$, fix a finite (\textit{i.e.} non-ideal) triangulation $t_i$ whose vertices are the points $M_k, N_k$ on $H_k, k=1, 2, ..., b$. Let $t_i$ be such that for every its edge $e'_j$ not in $\mathcal{E}_i$ ($i=1,2$), $f(e'_j)$ is entirely within one of the simple Euclidean polygons bounded by the images of edges of $\dot{F}$ from $\mathcal{E}_i$ under $f$.
For every $W_i$, there is more than one such triangulation; we choose $t_i$ so that the restriction of $f$ to every triangle is bijective, if such a triangulation exists (otherwise choose any triangulation).

   The vertices $v_1, ..., v_j$ of $\dot{F}$ are centers of the horospheres $H_1, ..., H_j$, whose sizes are defined by the corresponding diagram labels as explained in the previous section. For a vertex $v_i$ of $\dot{F}$ and two edges $e_j, e_{j+1}$ of $\dot{F}$ incident to $v_i$, let $F$ coincide with the hyperbolic plane defined by $e_j, e_{j+1}$ in the horoball neighborhood  of $v_i$ bounded by $H_i$. Consider a triangle $T$ from one of triangulations $t_i$. If $T$ is adjacent to a part of a face $F$ that lies inside a horoball (as defined in the previous paragraph), let $T$ be a part of $F$. If $T$ is not adjacent to any such part of any face, let $T$ be a part of any face to which it is adjacent along its other edges. This process defines $F$, which we call a ``straightened face" of $\Pi_1$.

  Note that there is a hyperbolic isometry that moves $v$ away from infinity, and places another vertex of $\Pi_1$ there. It can be used to choose the triangulation and therefore straightening of the remaining face of $\Pi_1$. Face identifications of $\Pi_1$ and $\Pi_2$ together with the choice of straightening of faces of $\Pi_1$ determine straightened faces of $\Pi_2$.  The collection of vertices, straightened edges and straightened faces of a polyhedron $\Pi_i$ will be denoted by $\partial\Pi_i$, and, when we want to exclude the vertex at infinity and adjacent edges and faces, by  $\partial\Pi_i-\{v\}$. Extend $f$ to be the vertical projection of   $\partial\Pi_i-\{v\}$ onto the plane $z=0$. When we write $\Pi_i$, $i=1,2$, throughout, we mean either an ideal hyperbolic polyhedron bounded by $\partial\Pi_i$, or a corresponding truncated polyhedron with cross-sections lying on the horospheres centered at the ideal vertices. Since the meaning will be clear from the context, we will use the same notation in both cases.

Below we formulate the conditions on the diagram labels that guarantee that the cusped polyhedra are properly embedded in $\mathbb{H}^3$. The conditions might seem unwieldy, but have a simple geometric meaning. Here and further we assume that the diagram labels satisfy the hyperbolicity equations.

\vv

a) Consider labels on two sides of an overpass of a crossing. There are four such labels, $u, v, u+1, v+1$. Moving in the
direction of the link orientation along the overpass, suppose the labels $u, u + 1$
come last, as in Fig.2. Then $\operatorname{Im}u>0$ holds iff $u$ is on the right with respect to our travel (Fig.2, left), and $\operatorname{Im}u<0$ holds if and only if it is on the left (Fig.2, right). This condition only needs to be satisfied for one arbitrarily chosen edge label of $D$ that is not purely real (if all labels are real, the condition is not satisfied automatically).

\begin{center}

\includegraphics[scale=0.68]{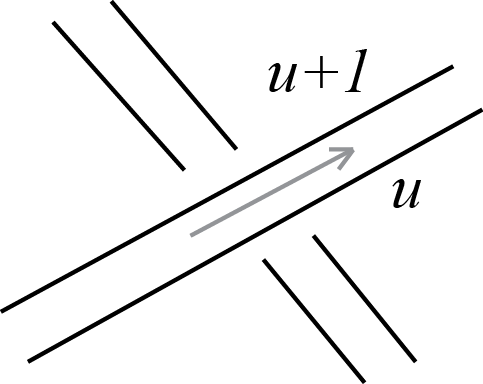} \hspace{0.5in}
 \includegraphics[scale=0.68]{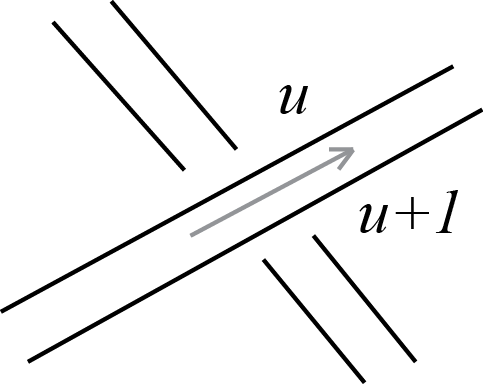} 

Fig.2
\end{center}

b) For every overpass/underpass that has two consecutive edge labels $u, v$ (Fig.3, left), and is not incident to a bigon, either $\operatorname{Im}(-(v+1)/u)>0$ and $\operatorname{Im}(-(u+1)/v)>0$ hold if $\operatorname{Im}u>0$, or $\operatorname{Im}(-v/(u+1))>0$ and $\operatorname{Im}(-u/(v+1)>0$ hold if $\operatorname{Im}u<0$.

  \begin{center}
  
  \includegraphics[scale=0.68]{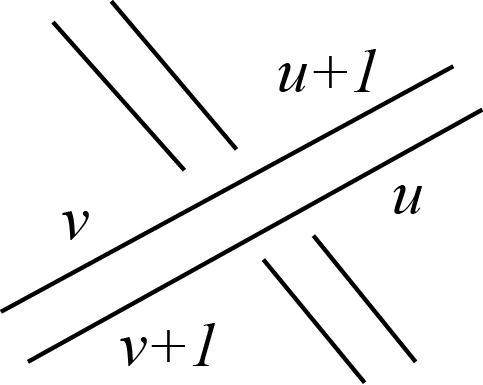} \hspace{0.5in}
   \includegraphics[scale=0.29]{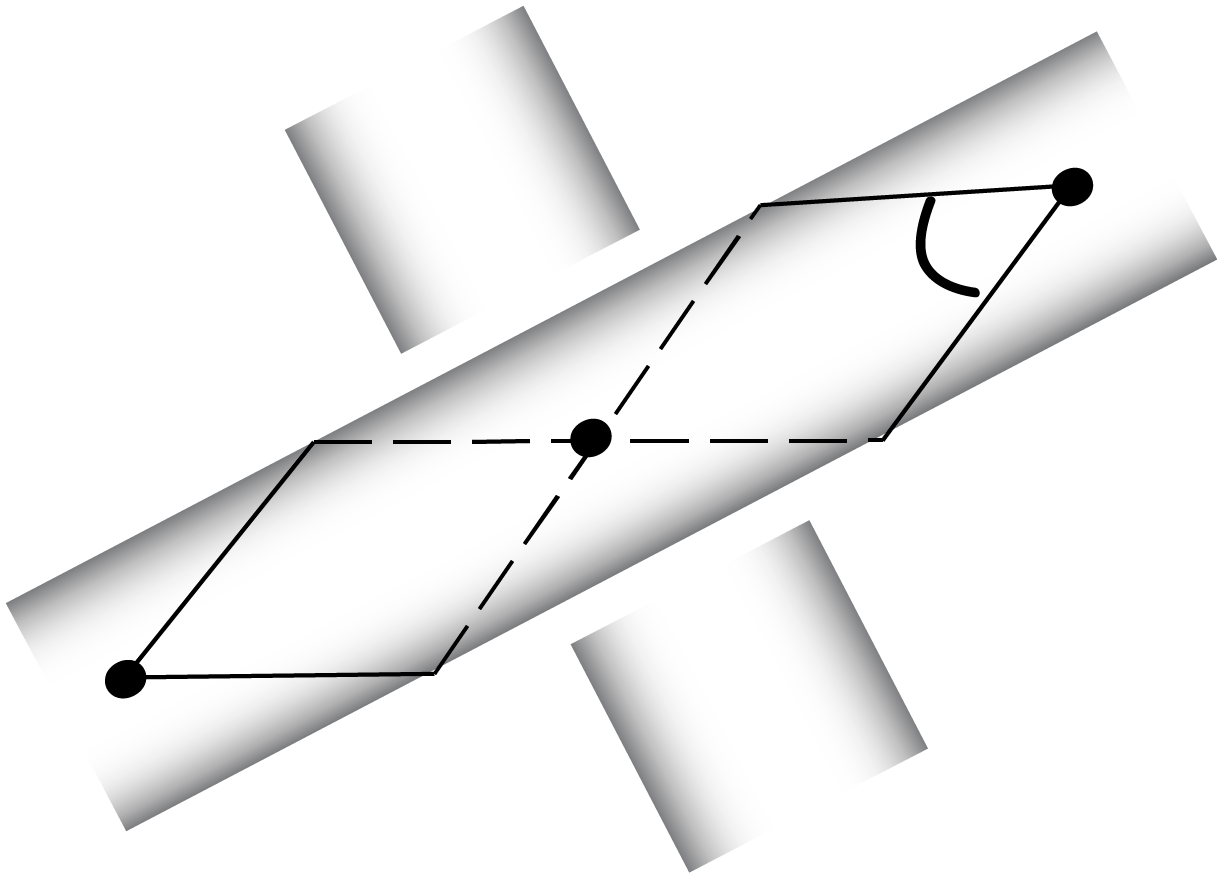} \hspace{0.5in}
  \includegraphics[scale=0.68]{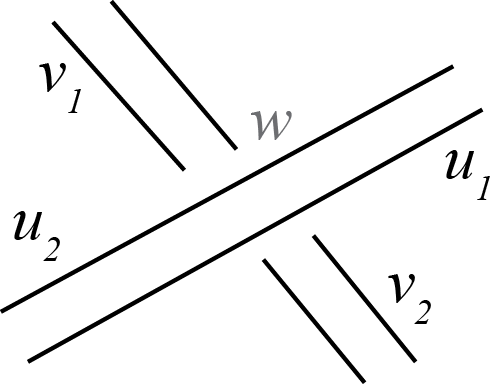}

  Fig.3
  \end{center}

c) There exists a crossing of $D$ for which the following holds. Suppose $w$ is the crossing label assigned to the crossing, and $u_i, i=1,2,$ are two edge labels on the adjacent overpass, and $v_j, j=1,2$ are two edge labels on the adjacent underpass (Fig. 4, right). Then at least one of the fractions $\frac{w}{u_i v_j}$, $\frac{w}{u_i(v_j+1)}$, $\frac{w}{v_j(u_i+1)}$ is not purely real for some $i,j$.

\vv

We will now explain the nature of these conditions.

\vv

\textbf{Lemma 3.1.} Conditions (a) and (b) imply the polyhedra $\Pi_1$ and $\Pi_2$ are cross-sectionally convex. Condition (c) guarantees that at least four horospheres are not all in the same plane, and not all crossing labels are 0.

\vv

\textbf{Proof.} In the polyhedron $\Pi_i, i=1,2$, there are at most four edges of $\mathcal{E}_i$ incident to an ideal vertex. The cross-section determined by these edges is often a Euclidean quadrilateral. The exception is a cross-section at an ideal vertex that has at least one incident polyhedral edge resulting from a bigon of $D$. Such a cross-section might be a triangle, or even a bigon that degenerates eventually into just one edge.

  Consider a quadrilateral cross-section. One may see one such cross-sectional tile of the torus boundary on Fig.3, middle, where a crossing of the thickened link is depicted  (two of the vertices of the tile are glued together underneath the overpass).

 If we look at the Fig.3, left, we can write the expressions for the angles of the cross-section in terms of the diagram labels. Two opposite angles of the cross-section are $\operatorname{arg}{\frac{u}{u+1}}$ and $\operatorname{arg}{\frac{-v}{-(v+1)}}$ (the corresponding quadrilateral is depicted in Fig.4, left). The condition (a) allows to choose the solution to the hyperbolicity equations that is consistent with the orientation conventions used in the definitions of the diagram labels (when choosing out of two Galois conjugates). The consistency for the other labels follows automatically from the equations after the correct choice of the solution is made. Note that if $u=a+bi$ and $\operatorname{Im}u=b>0$, then $\operatorname{Im} \frac{u}{u+1}=\frac{b}{(a+1)^2+b^2}>0$.  Therefore, the condition (a) also guarantees that $\operatorname{Im}\frac{u}{u+1}$ and $\operatorname{Im}\frac{v}{v+1}$ are positive, and therefore $\operatorname{arg}\frac{u}{u+1}$ and  $\operatorname{arg}\frac{v}{v+1}$ (corresponding to two opposite angles of the quadrilateral cross-section) are between $0$ and $\pi$.  In the case $\operatorname{Im}u=b<0$, $\operatorname{Im} \frac{u+1}{u}=\operatorname{Im}\frac{-b}{(a+1)^2+b^2}>0$, and the same angles are between 0 and $\pi$ again. The condition (b) similarly guarantees that two other opposite angles of the cross-section are between 0 and $\pi$ as well. 
  
 Consider a triangular cross-section (Fig.4, middle), resulted from collapsing a bigon region of the link diagram. If $\operatorname{Im}u=b>0$, then the interior angles correspond to $\operatorname{arg}\frac{u}{u+1},  \operatorname{arg}(u+1),  \operatorname{arg}-\frac{1}{u}$. The arguments of $\frac{u}{u+1}, u+1, -\frac{1}{u}$ are between 0 and $\pi$ (since $\operatorname{Im}\bigg(-\frac{1}{u}\bigg)=-b>0$).  If $\operatorname{Im}u=b<0$, then the interior angles correspond to $\operatorname{arg}\frac{u+1}{u},  \operatorname{arg}(-u),  \operatorname{arg}\frac{1}{u+1}$, and are   between 0 and $\pi$ again. Therefore, the conditions (a) and (b) together imply cross-sectional convexity of the two polyhedra.
  
  \begin{center}
  \includegraphics[scale=0.53]{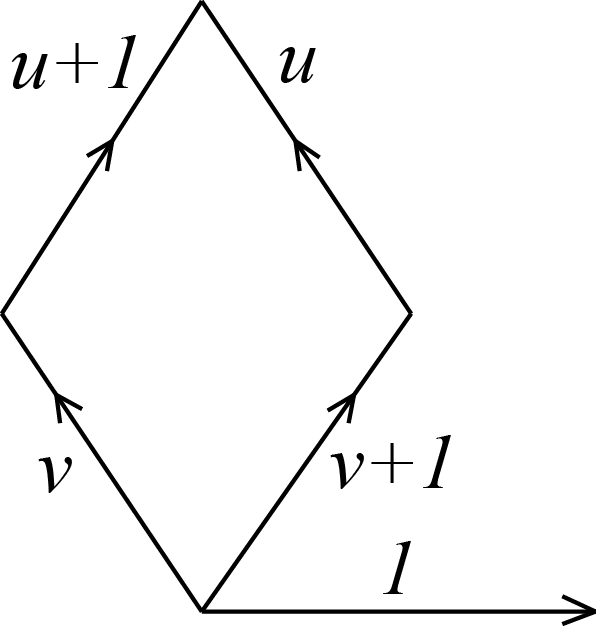} \hspace{0.5in}
  \includegraphics[scale=0.72]{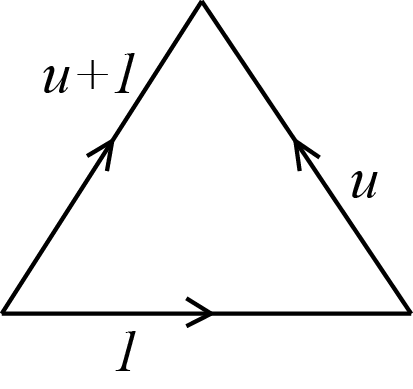}
  \hspace{0.5in}
  \includegraphics[scale=0.64]{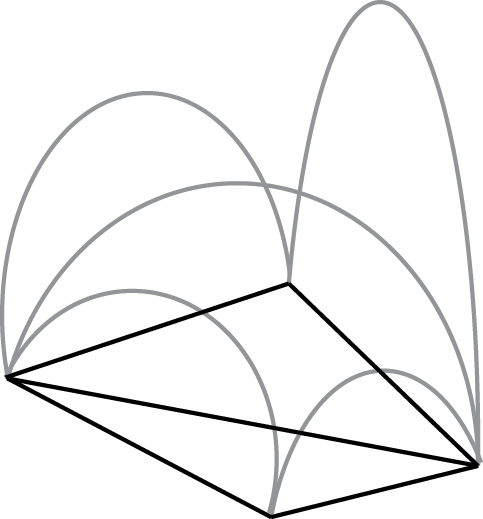}

   Fig.4 
  \end{center} 

 The condition (c) guarantees that we have at least four horospheres that are not all located in the same plane, since the quotient $\frac{w}{u_i v_j}$ is the cross-ratio of the centers of the corresponding horospheres (possibly, with a negative sign). For the details and proof, see Sections 2-4 in \cite{method}. It also ensures that not all crossing labels are 0, which would force the system of hyperbolicity equations to be degenerate, of the form $0=0$. $\Box$
\vv

We can now proceed to the main result of this section.

\vv

\textbf{Proposition 3.2.} Suppose that for an alternating link $L$ there exists a solution to the hyperbolicity relations such that the conditions (a)--(c) are satisfied. Then each of ideal polyhedra $\Pi_1, \Pi_2$ is properly embedded in $\mathbb{H}^3$, \textit{i.e.} any two straightened faces either are identified, or are disjoint, or intersect in a connected sequence of edges (including the vertices that are the edges' endpoints), or in an ideal vertex only.

\vv

The rest of this section is devoted to the proof of the proposition. The proof consists of a number of observations about the geometric nature of the polyhedra.
 
Note that it is enough to prove the statement solely for the faces of $\Pi_1$, since the construction of $\Pi_2$ is similar to that of $\Pi_1$. Once the Proposition 3.2 holds for any two faces of $\Pi_1$ and any two faces of $\Pi_2$, it holds for any pair of faces where one face is from $\Pi_1$ and the other face is from $\Pi_2$ by the construction as well. 

By Lemma 3.1, $\Pi_1$ is cross-sectionally convex. This implies, in particular, that for any horosphere $H$ centered at an ideal vertex $v$ of $\Pi_1$, the interior angles (not exterior) of the cross-section of $\Pi_1$ on $H$ correspond to the dihedral angles of $\Pi_1$ in the horoball neighborhood of $v$ bounded by $H$. It also implies that for any edge $e$ of $Q$, the intersection of $H$ and the interior of $\Pi_1$ is on one side of $e$ (not on both), and is inside $Q$. Note also that the boundary of $\Pi_1$ is simply connected, and $\partial\Pi_1-\{v\}$ is connected.

\vv

\textbf{Lemma 3.3.} Under the hypothesis of Proposition 3.2, there is no point $p$ in $f(\partial\Pi_1-\{v\})$ such that $f^{-1}(p)$ consists of at least two distinct points.  

\vv

\textbf{Proof of Lemma 3.3.}  Suppose the contrary.

Assume $p$ is in the horoball neighborhood $H$ of an ideal vertex of $\partial\Pi_1-\{v\}$. Our definition of straightened faces implies that $f$ is not bijective in such a horoball neighborhood in two cases. Firstly, if a face of $\Pi_1$ inside $H$ lies in a vertical plane. But then this face of $\partial\Pi_1-\{v\}$  has $\infty$ among its vertices, a contradiction. Secondly, if there are two faces of $\Pi_1-\{v\}$ inside $H$, one directly above another. Since $\partial\Pi_1-\{v\}$ is connected, this leads to two layers of faces, one above another, outside the horospherical neighborhood of the ideal vertex.

 Now assume $p$ is outside the horoball neighborhoods of ideal vertices of $\partial\Pi_1-\{v\}$.  Either two distinct points of $f^{-1}(p)$ belong to one triangle that is a part of a straightened face of $\Pi_1$, or to two distinct triangles.

First, suppose distinct points from $f^{-1}(p)$ belong to the same triangle $T$ of a face $F$. Since $f$ is a vertical projection, $T$ is vertical as well (\textit{i.e.} $T$ is perpendicular to the plane $z=0$). Recall that for the collection of vertices and straightened edges $\dot{F}$ corresponding to the vertices and edges of a face of $\bar{\Pi}_1$, the triangulation is chosen so that every triangle projects bijectively under $f$, if such a triangulation exists. Such a triangulation does not exist only if there are two distinct consecutive edges of $\dot{F}$ that lie in one vertical plane. Denote them by $e_1, e_2$.

The edges $e_1, e_2$ are incident to an ideal vertex, say $v$, and there is a horosphere centered at $v$, say $H$. Denote the cusp cross-section of $\Pi_1$ on $H$ by $Q$. Two vertices of $Q$ that are on $e_1, e_2$ lie in the vertical plane. The other two vertices of $Q$ must be outside the plane and on the same side from it to satisfy cross-sectional convexity. This implies that there are two levels of  $\partial\Pi_1-\{v\}$, one above another, on one side of $T$.

If two distinct points of $ f  ^{-1}(p)$ do not belong to the same triangle, they belong to two different triangles, and there are two levels of $\partial\Pi_1-\{v\}$, one above another, as well.

 The assumption that there are two levels of faces of $\partial\Pi_1-\{v\}$ leads to the following cases.

\v

 Case 1.  A face of $\Pi_1$ incident to $\infty$ meats a face of the lower level,  say $K_3$, and a face of the upper level, say $K_4$ simultaneously at a certain vertex. We will denote the horosphere centered at this vertex by $H_{34}$. This is depicted on Fig.5, left, with the faces of the upper level in black, and the faces of the lower level in grey (the triangulations of the faces are not pictured). 
  
  \begin{center}
  
  \includegraphics[scale=0.85]{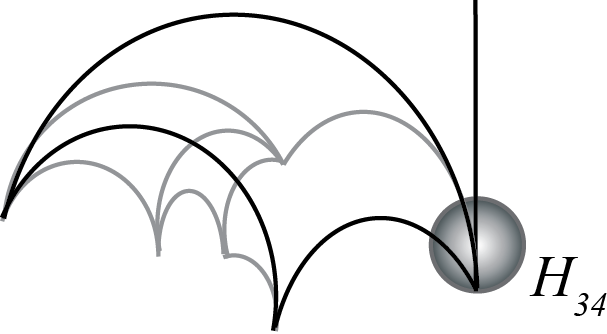} \hspace{0.6in}
  \includegraphics[scale=0.45]{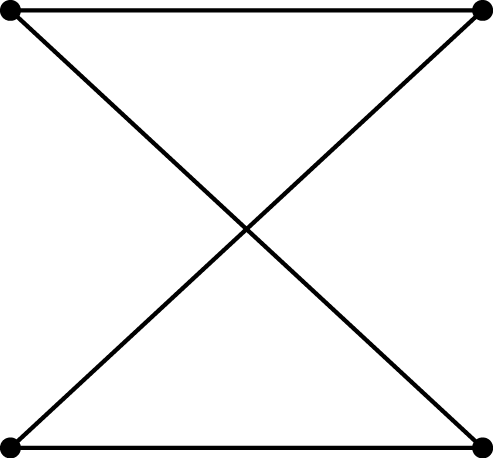}

  Fig.5 
  \end{center}

Two edges of $\mathcal{E}_i-{v_i}$ incident both to the center of $H_{34}$ and to $K_3$ determine a plane that intersects $H_{34}$ in an edge $e$. The vertices of the cross-section on $H_{34}$ that are not the endpoints of $e$ are on the opposite sides of the Euclidean line defined by $e$ on $H_{34}$ (Fig.5, right). This contradicts cross-sectional convexity.

Note that the same argument leads to a contradiction if at a vertex of $\partial\Pi_i-\{v\}$, faces of three or more levels meet. The same argument also gives a contradiction if two layers of faces meet at a edge $e_p$ of $\mathcal{E}_i$, and the faces are adjacent to a single layer of faces at $e_p$.
 
 \v
 
 Case 2. Faces meeting $\infty$ are incident at their vertices only to faces on the lower level. Denote a face with a vertex at $\infty$ by $K_4$.
 Then we arrive at one of the two scenarios. Either $K_4$ has interior of the polyhedron on both sides near the plane $z=0$ (as on Fig.7, left). Then on one of the horospheres adjacent to $K_4$ and not located at $\infty$, the cross-section looks like Fig.5, right, again. A contradiction. Or, in the second scenario, none of the faces that are incident to $\infty$ and to lower level faces simultaneously have interior of the polyhedron on both sides (as on Fig.7, right). Then the exterior angles of the cross-section correspond to the dihedral angles of the polyhedron in the corresponding horoball neighborhood of the vertex. A contradiction.

 \v
 
 Case 3. Faces meeting $\infty$ meet only faces in the upper level at vertices. However there is
 a lower level of faces elsewhere. We then arrive at one of two scenarios. The first scenario is that the faces of upper and lower level meet at some edges of $\partial\Pi_i-\{v\}$ and are adjacent to a single layer of faces through those edges (and this single layer is then adjacent a face with a vertex at $\infty$). We then use the argument from case 1. The second scenario is to have a face with a vertex at $\infty$ adjacent to the upper level of faces which then becomes the lower level. Then there is an ideal vertex $q$ that is a vertex of the face $K_5$ at least a part of which is in the upper level, and of the face $K_6$ at least a part of which is in the lower level. Denote the horosphere centered at $q$ by $H_{56}$. $K_5$ intersects $H_{56}$ in an edge $e$. Then the interior of the polyhedron is on both sides of $e$, since it is above the faces of the upper level adjacent to $\infty$, and right above $K_6$ as well. A contradiction.
 
 \v

Case 4. Faces meeting $\infty$ meet a single layer of faces. Recall that the boundary of the polyhedron is connected, and we have two layers of faces elsewhere. Hence the faces of upper and lower levels either meet at some edges of $\partial\Pi_i-\{v\}$ and are adjacent to a single layer of faces at those edges, or intersect and switch levels. In the former situation, the argument from case 1 leads to a contradiction. In the latter situation, the argument from the second scenario of case 3 leads to a contradiction. $\Box$

\vv

 \textbf{Corollary 3.4.} The image of the edges of $\mathcal{E}_i$, $i=1,2$ that bound a face of $\Pi_i$  under $f$ is a simple polygon. Moreover, a restriction of $f$ to every straightened  face of $\Pi_i$ is a bijection. 
 
 \vv

\textbf{Proof of Proposition 3.2.} If there are two distinct faces of $\partial\Pi_1-\{v\}$ with intersecting interiors, then there is a point $p$ in the image of one of the faces under $f$ such that its preimage consists of at least two distinct point. A contradiction to the claim.

  \begin{center}

  \includegraphics[scale=0.5]{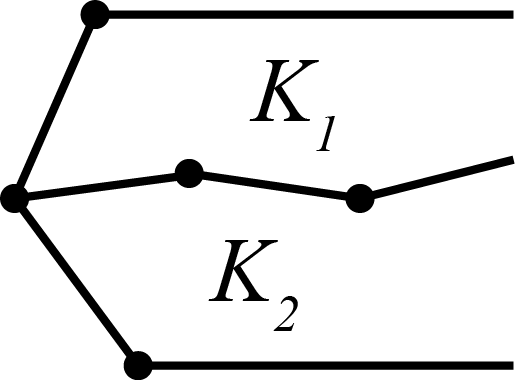} \hspace{0.6in}
  \includegraphics[scale=0.35]{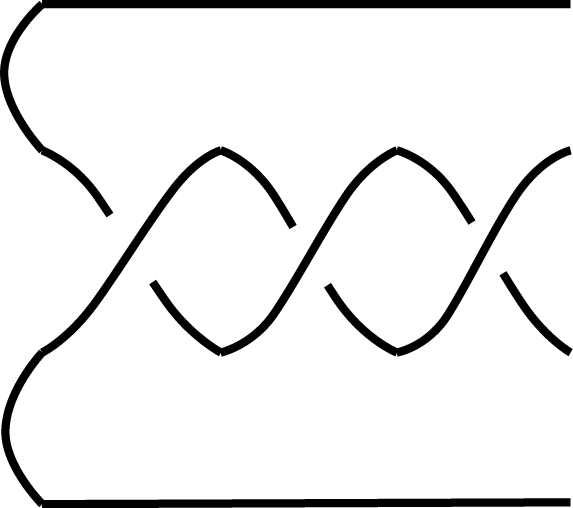}
  
  Fig.6 
  \end{center}

Suppose there are two distinct faces $K_1, K_2$ of $\Pi_1$, whose boundary intersects in more than just a vertex or a sequence of connected edges (as on Fig.6, left). $K_1, K_2$ result from a reduced and twist-reduced alternating diagram. Then Menasco's reduced alternating link diagram construction of the polyhedra implies that then there are two faces of $D$ sharing more than a connected sequence of edges or a crossing. In a reduced and twisted-reduced alternating diagram, this is possible only if two adjacent non-bigon regions have a number of bigons between them, as in Fig.6, right, and share several crossings. But bigons are collapsed in Menasco's construction, and result in a row of consecutive edges in the corresponding polyhedra, that the faces share. A contradiction. $\Box$

\begin{center}

\includegraphics[scale=0.85]{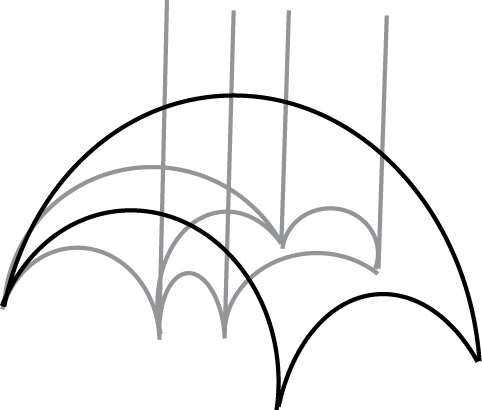} \hspace{0.5in}
\includegraphics[scale=0.85]{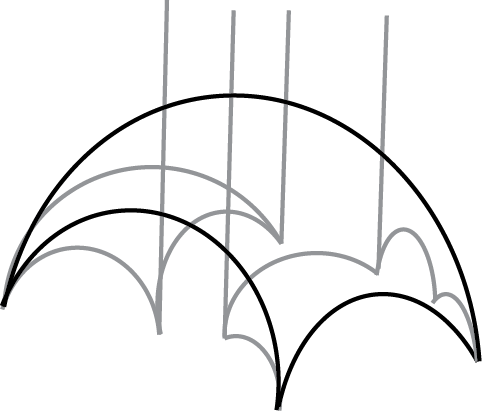}

Fig.7 
\end{center}

\section{Ideal partially flat geodesic triangulations}

Thurston provided sufficient conditions for an ideal triangulation of a finite volume 3-manifold to be geodesic, which can be expressed in gluing and completeness equations (\cite{Thurston}) together with a certain restriction on a solution. In particular, all the tetrahedra should have positive volume. Petronio and Weeks showed that this can also be achieved if an ideal triangulation is partially flat and satisfies the completeness and consistency conditions (\cite{PetronioWeeks}). 

In this section, we will show that any cross-sectionally convex polyhedra in $\mathbb{H}^3$ can be subdivided into a partially flat geodesic triangulation. This implies that the conditions (a)--(c) guarantee the existence of the complete hyperbolic structure of $S^3-L$, and the correspondence described in Section 2 induces a developing map.

Let us recall the nature of the conditions that make a triangulation an ideal geodesic partially flat triangulation. 
Consider a Euclidean cross-section of an ideal tetrahedron. Suppose that the (complex)
translations corresponding to the sides of the cross-section are $u, v, -(u + v)$ (Fig.8). Any of \h $-\frac{v}{u}$, \h $\frac{u}{u+v}$, \h $\frac{u+v}{v}$ \h can be taken as the shape $z$ of the tetrahedron (sometimes also called tetrahedral parameter or modulus), and their arguments correspond to the interior Euclidean angles of the cross-section, as well as to the interior dihedral angles of the tetrahedron.

\vv

\textbf{Definition 4.1.} The completeness and consistency conditions for an ideal triangulation $\tau$ are as follows.

\begin{center}

\includegraphics[scale=0.45]{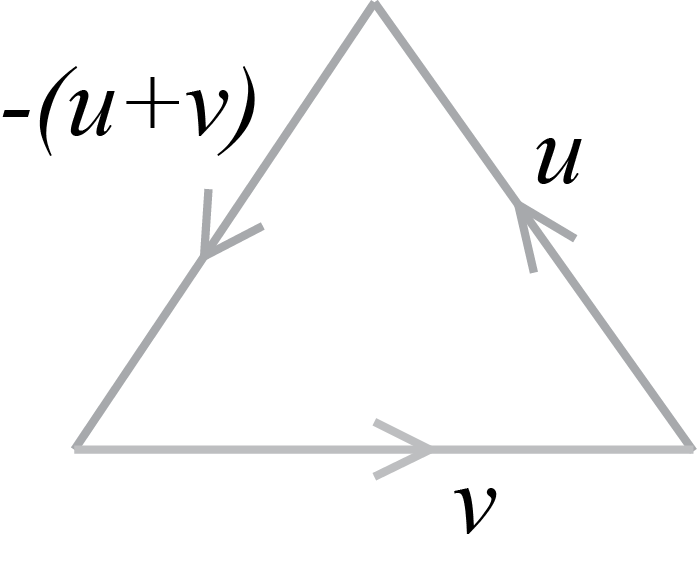} 

Fig.8

\end{center}

1) In every tetrahedron, three pairs of opposite edges correspond to $z$, $1-\frac{1}{z}$, $\frac{1}{1-z}$ in the way indicated above.

2) At every edge $e$ of $\tau$, after all the faces are identified in pairs, the shapes of the tetrahedra glued at $e$ satisfy $z_1 z_2 ... z_k = 1$.

3) For every tetrahedron, its shape $z$ satisfies $\operatorname{Im}z\geq0$, and not all $z_i, i=1,2...,n$, are 0.

4) The metric is complete, \textit{i.e.} cross-sectional triangles glued together at each ideal vertex must fit together to give a closed Euclidean surface

\vv
Note that the condition $\operatorname{arg}z_1+ \operatorname{arg}z_2+... +\operatorname{arg}z_k=2\pi$ is often included (as a part of consistency conditions), but one can prove that if both (2) and (4) are satisfied, it is redundant.

\vv

\textbf{Theorem 4.2.} \h Suppose that for an alternating link $L$ with a reduced alternating diagram $D$ edge and crossing labels satisfy the conditions (a)--(c). Then there exists a partially flat geodesic triangulation $\tau$ that induces the complete hyperbolic structure on $S^3-L$. 

\vv

Most of this section is devoted to the proof of Theorem 4.2. 

Under the hypothesis of the theorem, two straightened polyhedra $\Pi_1, \Pi_2$ that correspond to Menasco's decomposition of $S^3-L$ are properly embedded in $\mathbb{H}^3$ by Prop 3.1. First, we construct a triangulation of $\Pi_1 \cup \Pi_2$ that is properly embedded in $\mathbb{H}^3$. Recall that $f$ is the vertical projection described in the previous section. To begin, take the collection of ideal vertices and straightened edges $\dot{F}$ corresponding to a face $F$ of $\Pi_i-\{v\}$. Subdivide $\dot{F}$  into (ideal) triangles using the existing ideal vertices of $\dot{F}$ only. Additionally, do it so that for any new edge $e$ subdividing $\dot{F}$ of $\Pi_i-\{v\}, i=1,2$, $f(e)$ lies entirely in $f(F)$ (this can always be done as a consequence of Cor.3.3). Do this for every face of $\Pi_i-\{v\}, i=1,2$. Additionally, for every face $F$ of $\Pi_i-\{v\}, i=1,2$, incident to $\infty$, subdivide $\dot{F}$ by edges incident to $\infty$ (i.e. by vertical geodesics from other vertices of $\dot{F}$). Denote the resulting polyhedra with triangular faces by $\Pi'_1, \Pi'_2$. 

Recall that a collection of edges and vertices in $\mathbb{H}^3$ corresponding to a face $\bar{F}$ of $\bar{\Pi}_1$ was denoted by $\dot{F}$. As mentioned before, the ideal vertices of $\dot{F}$ do not necessarily lie in one hyperbolic plane. Therefore, after we subdivided every such collection $\dot{F}$ in a new way, we have to check that the resulting polyhedra $\Pi'_i, i=1,2,$ are properly embedded in $\mathbb{H}^3$ as well (\textit{i.e.} the new triangular faces do not ``cut" through each other). That is the content of the next lemma.

Fix $i$ and consider a cross-section $Q'$ of an ideal vertex of the polyhedron $\Pi_i'$ resulting from the subdivision. Suppose $Q'$ is situated on the horosphere $H$. $Q'$ corresponds to a 4-sided cross-section $Q$ of $\Pi_i$ with the vertices $A, B, C, D$ in the following way: $Q'$ has four vertices $A, B, C, D$ of $Q$, and several more vertices resulting from the subdivision of faces of $\Pi_i$. We will refer to these vertices as to the ``new" vertices of $Q'$.

\vv

\textbf{Lemma 4.3.} \h Suppose $E$ is a new vertex of the cross-section $Q'$. Suppose further $E$ resulted from the subdivision of the face of $\Pi_i$ that is incident to the edge $CD$ of $Q$. Then $E$ lies either in the interior of $Q$, or on $CD$, or in the part of $H$ bounded by the lines $CD, AD, BC$ and on the other side of the line $CD$ from $AB$ (the area where $E$ may lie is shown in grey on Fig.9).

\vv
\textbf{Proof of Lemma 4.3.} Suppose first $\dot{F}$ does not meet $v$. Since $f$ is a homeomorphism on faces of $\Pi_i-\{v\}$ (by Proposition 3.2), any new edge of $\Pi_i'$ resulting from the subdivision of $\dot{F}$ does not intersect the interior of any face of $\Pi_i$ besides $F$.

Now suppose $\dot{F}$ has $v$ among its vertices. Then the new edges resulting from the subdivision of $\dot{F}$ are all vertical geodesics. Since there are no two levels of faces anywhere (see the proof of Lemma 3.3), these new edges do not pierce the interior of any face (besides $F$) as well.

Denote the geodesics that pierce $H$ at the points $A, B, C, D, E$ by $\alpha, \beta, \gamma, \delta, \eta$ respectively. If $E$ lies outside the specified area, $\eta$ pierces a face adjacent to either $AD$ or $AB$ or $BC$, a contradiction. $\Box$

\begin{center}
\includegraphics[scale=0.72]{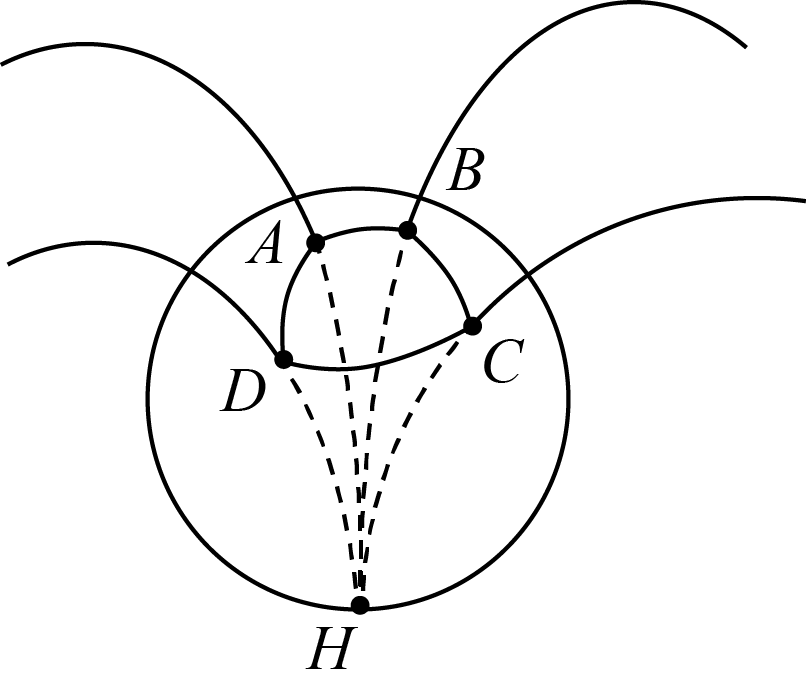} \hspace{1cm}
\includegraphics[scale=0.72]{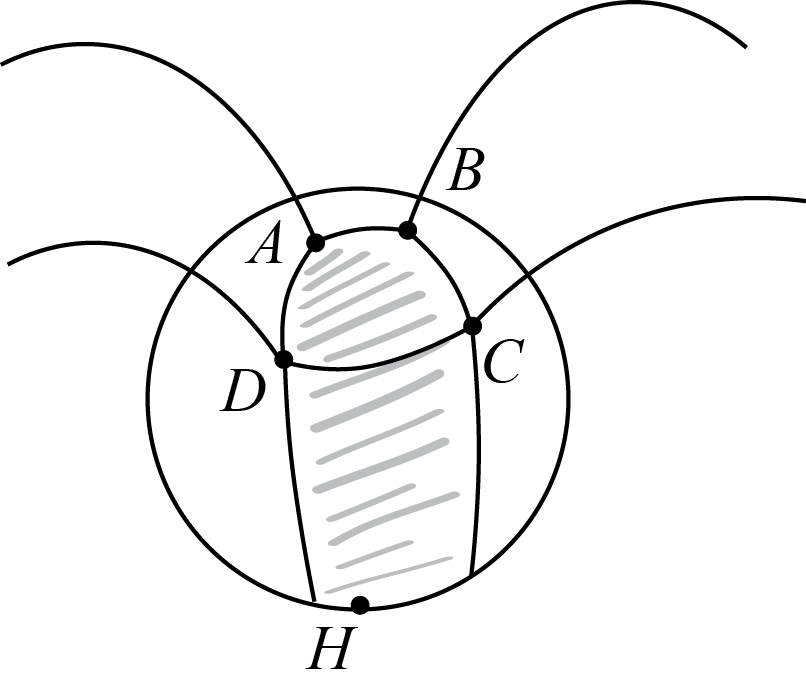}

Fig.9 
\end{center}

All the faces of $\Pi'_1, \Pi'_2$ are triangular faces, and the previous lemma implies that $\Pi'_1, \Pi'_2$ are properly embedded in $\mathbb{H}^3$. We will now construct an ideal partially flat triangulation of $\Pi'_1, \Pi'_2$, properly embedded in $\mathbb{H}^3$.

 Triangulate $\Pi'_1, \Pi'_2$ using the existing ideal vertices only, and so that the interiors of the triangular cross-sections of tetrahedra do not intersect. Denote such a triangulation by $\tau$. Note that such $\tau$ always exists and is not unique; an example is the triangulation suggested by Thurston (\cite{Thurston}) and Menasco (\cite{Menasco}). 
 
 \vv
 
 \textbf{Lemma 4.4.} Any tetrahedron in $\tau$ is either flat, or lies entirely inside the polyhedron $\Pi'_1$, $i=1,2$, or lies entirely outside.
\vv

 \textbf{Proof of Lemma 4.4.} The subdivision of cross-sections of $\Pi'_i$ corresponds to the subdivision of the polyhedron $\Pi'_i$ in the following way. Suppose we subdivided the cross-section $Q'$ with the consecutive vertices $A, B, C, D, E$ on the horosphere $H$ by a new edge $EC$. Suppose $P$ is the center of $H$. Suppose also that the vertices $A, B, C, D, E$ result from geodesics $\alpha, \beta, \gamma, \delta, \eta$ respectively (as before) piercing $H$, and that these geodesics are edges of $\Pi_i'$. Then the plane $P_{EC}$ defined by the geodesics $\gamma, \eta$ intersects $Q'$ in $EC$. Suppose also that $P_1, P_2$ are the ideal vertices that are adjacent to $\gamma, \eta$ respectively on the sides opposite of $P$. We subdivide the cross-section $Q'$ of $\Pi'_i$ by an edge $EC$ if we also subdivide $\Pi'_i$ by a new triangular face that lies in the plane $P_{EC}$ and has vertices $P, P_1, P_2$ (possibly with a new edge between $P_1$ and $P_2$).

After the triangulation process, all cross-sections of $\Pi'_i$ are subdivided into triangles. If there is a tetrahedron $T$ of $\tau$ that lies partially outside and partially in the polyhedron $\Pi'_i$ in $\mathbb{H}^3$, then at least one of the faces (say, a face $F$) of $T$ is partially outside and partially in $\Pi'_i$. Then in the corresponding cross-section we have a subdividing edge that is partially inside and partially outside the cross-section, \textit{i.e.} the subdividing edge intersects edges of the cross-section more than just in their endpoints. This is not possible, since by Lemma 4.3 $Q'$ is a simple polygon, and since by our construction the triangulation of $\Pi_i'$ corresponds to triangulating this simple polygon. Similar argument shows that the interior of any truncated tetrahedron is not intersected by another truncated tetrahedron. $\Box$

\vv

\textbf{Proof of Theorem 4.2.} Lemma 4.4 implies that $\tau$ is an ideal partially flat triangulation. Additionally, the condition (c) implies that not all tetrahedra have 0 volume (as explained in Lemma 3.1). Hence, the condition (3) from Definition 4.1 is satisfied. Now let us check the conditions (1), (2) and (4).

Let us look at the condition (2) of Definition 4.1. Choose an edge $e$ of $\tau$, and let $T_1, ...., T_k$ be tetrahedra glued at $e$. Their shapes are the ratios of the edge labels $z_1=u_1/u_2$, $z_2=u_2/u_3$, ..., $z_{k-1}=u_{k-1}/u_k$, $z_k=u_k/u_1$, where each $u_i$ connects the point where $e$ pierces $H_1$ with the points where other edges of $T_1, T_2, ... , T_k$ pierce $H_1$ (Fig.10, left). None of $u_i$ is 0 due to the conditions (a)-(c), and the product of such shapes then satisfies $z_1 .... z_k =1$.

Turn our attention to the condition (1) of Definition 4.1. Consider a tetrahedron $T$ in the triangulation $\tau$, with one of its ideal vertices being the center of a horosphere $H_1$. Suppose the (complex)
translations corresponding to the sides of a cross-section lying at the horosphere $H_1$ are $u_1, v_1, -(u_1 + v_1)$. Let the geodesic edges of $T$ be denoted by $\gamma_i, i=1,2, ..., 5$ as on Fig.10, right. Then the shapes associated with the geodesics $\gamma_1, \gamma_2$ and $\gamma_3$ are \h $-\frac{v_1}{u_1}$, \h $\frac{u_1}{u_1+v_1}$, \h $\frac{u+v}{v}$. Denote the first shape by $z$, and then the other two shapes are $1-\frac{1}{z}, \h \frac{1}{1-z}$. The corresponding Euclidean angles of the cross-section are arguments of the shapes, and so are three corresponding dihedral angles. 

We need to check whether opposite dihedral angles in $T$ agree. For this, it is enough to check that the shapes agree. Suppose, another cross-section of $T$ lying on the horosphere $H_2$ has sides $u_2, v_2, -(u_2 + v_2)$.

Diagram labels satisfy the region relations, and the relations are obtained by composing the isometries rotating (truncated) hyperbolic polygons. Each polygon is a preimage of the boundary of a region of $D$. We may use faces of a triangulation instead of the regions (the labels are defined so that they satisfy the relations coming from rotating these faces as well). All these faces are three-sided, which makes the relations particularly simple. 

\begin{center}

\includegraphics[scale=0.66]{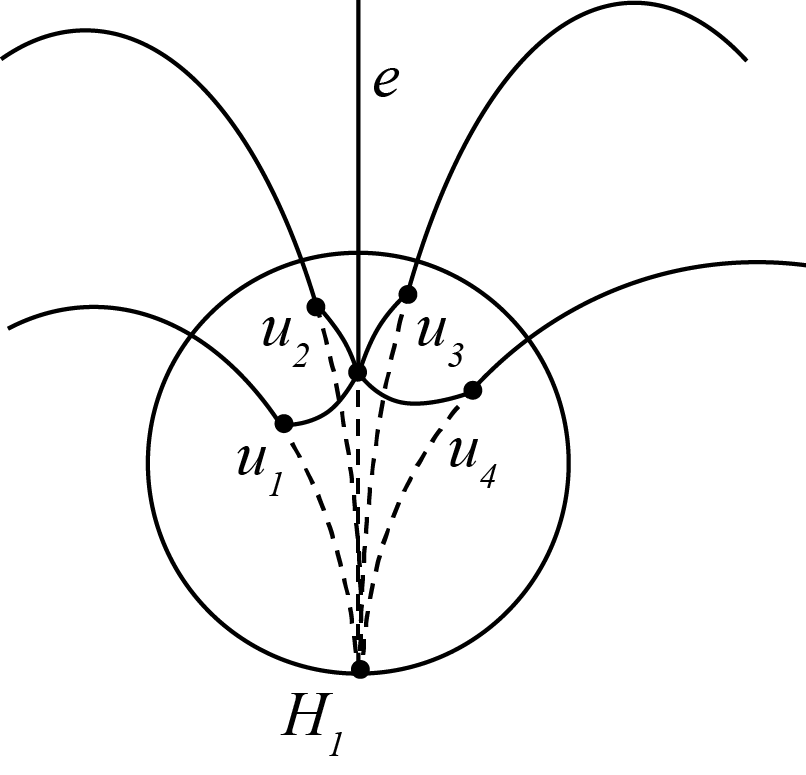} \hspace{0.3 in}
\includegraphics[scale=0.68]{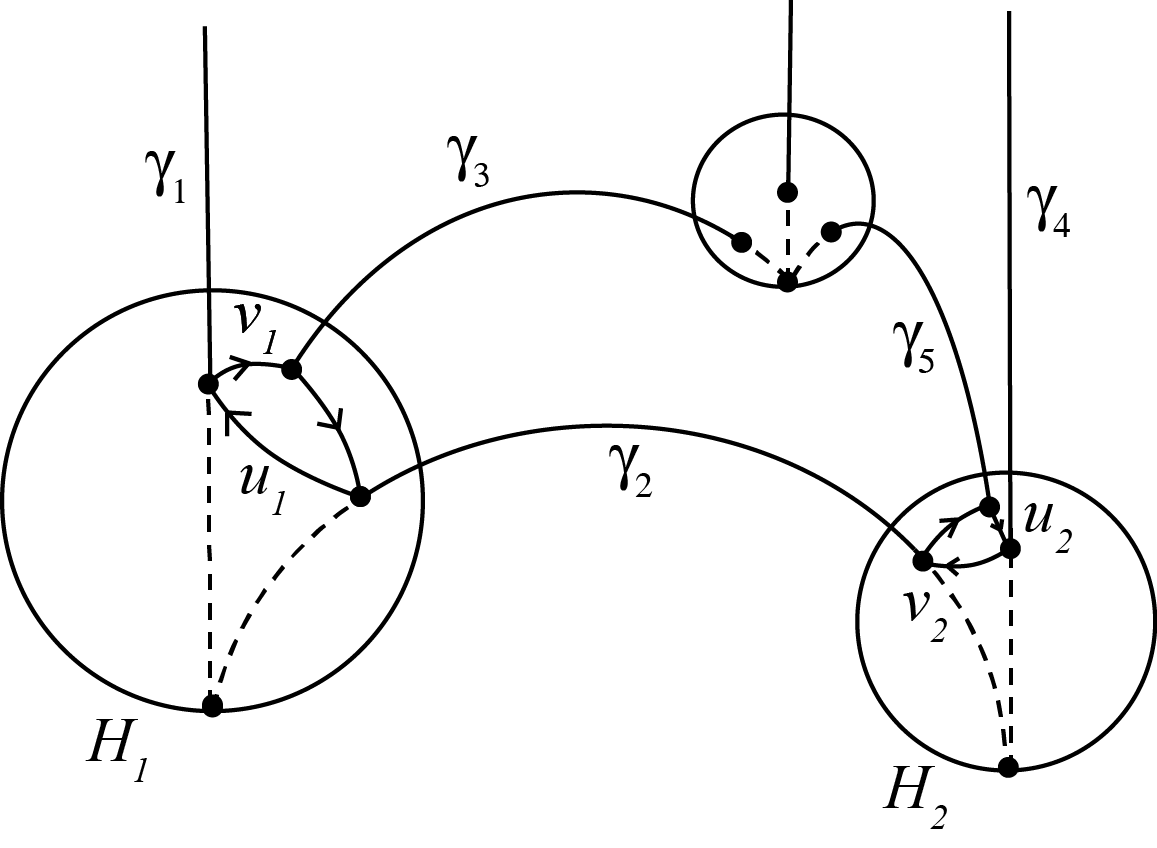} 
 
Fig.10

\end{center}

Let us show that the shape of $T$ associated with $\gamma_3$ is equal to the shape of $T$ associated with $\gamma_4$. The former is $\frac{u_1+v_1}{v_1}$, while the latter is $\frac{-v_2}{u_2}$. 

From the 3-sided polygon that is a face of $T$ determined by geodesics $\gamma_1, \gamma_2, \gamma_4$: $\frac{w_2}{u_1v_1}=1$, and therefore $v_2=\frac{w_2}{u_1}$. From the 3-sided polygon that is a face of $T$ determined by $\gamma_3, \gamma_2, \gamma_5$: $\frac{w_2}{(u_1+v_1)(u_2+v_2)}=1$, and therefore $u_2+v_2=\frac{w_2}{u_1+v_1}$. Then, from the cross-section of $T$ on the horosphere $H_2$, $u_2=(u_2+v_2)-v_2=\frac{w_2}{u_1+v_1}-\frac{w_2}{u_1}$. 

Substitute the shape of $T$ associated with $\gamma_4$ into the above expressions for $v_2$ and $u_2$. After routine simplifications, we obtain $\frac{-v_2}{u_2}=-\frac{u_1+v_1}{v_1}$, which is exactly the shape of $T$ associated with $\gamma_3$. Similarly one can show that other pairs of shapes of $T$ agree.

Lastly, let us check the condition (4) of Definition 4.1, namely that Euclidean cross-sections of all tetrahedra incident to a particular ideal vertex $v_1$ of $\Pi_i, i=1,2$, form a closed Euclidean surface. The tetrahedra are glued at an ideal edge incident to $v_1$, say $e$. Denote the other ideal vertex of $e$ by $v_2$. Denote the horopshere centered at $v_1$ by $H_1$. 

Recall that every ideal vertex of $\tau$ in $\mathbb{H}^3$ corresponds to an overpass or an underpass of the diagram $D$. Every edge connecting the vertices of $\tau$ lies on a geodesic that connects the centers of the corresponding horospheres in $\mathbb{H}^3$. Hence, there is a quadrilateral cross-section of $\Pi_i$ on $H_1$ resulting from connecting $H_1$ with four neighboring horospheres by four distinct arcs (Fig.11, left). Let $u_1, u_2, u_3,
u_4$ be the corresponding edge labels as in the figure. Since $u_1, u_2, u_3, u_4$ satisfy the definition of edge labels, they are determined by the Euclidean translations between the points where the corresponding arcs pierce $H_1$. 

Suppose first $e$ is an edge of $\Pi_i, i=1,2$. Then there is a vertex $Q$ of the quadrilateral cross-section on $H_1$ such that $Q$ lies on $e$. Three such cross-sections meet at $Q$ as in the Fig.3, center. The edges of these cross-sections are $u_1, u_2, u_3, u_4$ by the definition of edge labels and our construction of the polyhedra $\Pi_i, i=1,2$. Therefore, the cross-sections meet to form a closed Euclidean surface.

We triangulated $\Pi_1\cup\Pi_2$ by adding more arcs (edges) between the existing ideal vertices. For an ideal vertex $v_1$, all the arcs of the triangulation start at $H_1$ and end at one of the neighboring horospheres, which correspond to other overpasses and underpasses. Suppose now $e$ is not an edge of $\Pi_i, i=1,2$, and was added in the triangulation process. Define $u'_3$ and $u''_3$ to be complex numbers corresponding to the translations between the point where $e$ pierces $H_1$, and the points where the arcs immediately to the right and immediately to the left to $e$ pierce $H_1$. The numbers $u'_3$ and $u''_3$ are uniquely determined by the previous labels and the choice of the ideal vertices (and hence location on the plane $z=0$) for the incident tetrahedra. In particular, $u'_3+u''_3=u_3$ as in Fig.11, right. Hence, triangular cross-sections of all tetrahedra fit together to form a closed Euclidean surface. 

\begin{center}

\includegraphics[scale=0.45]{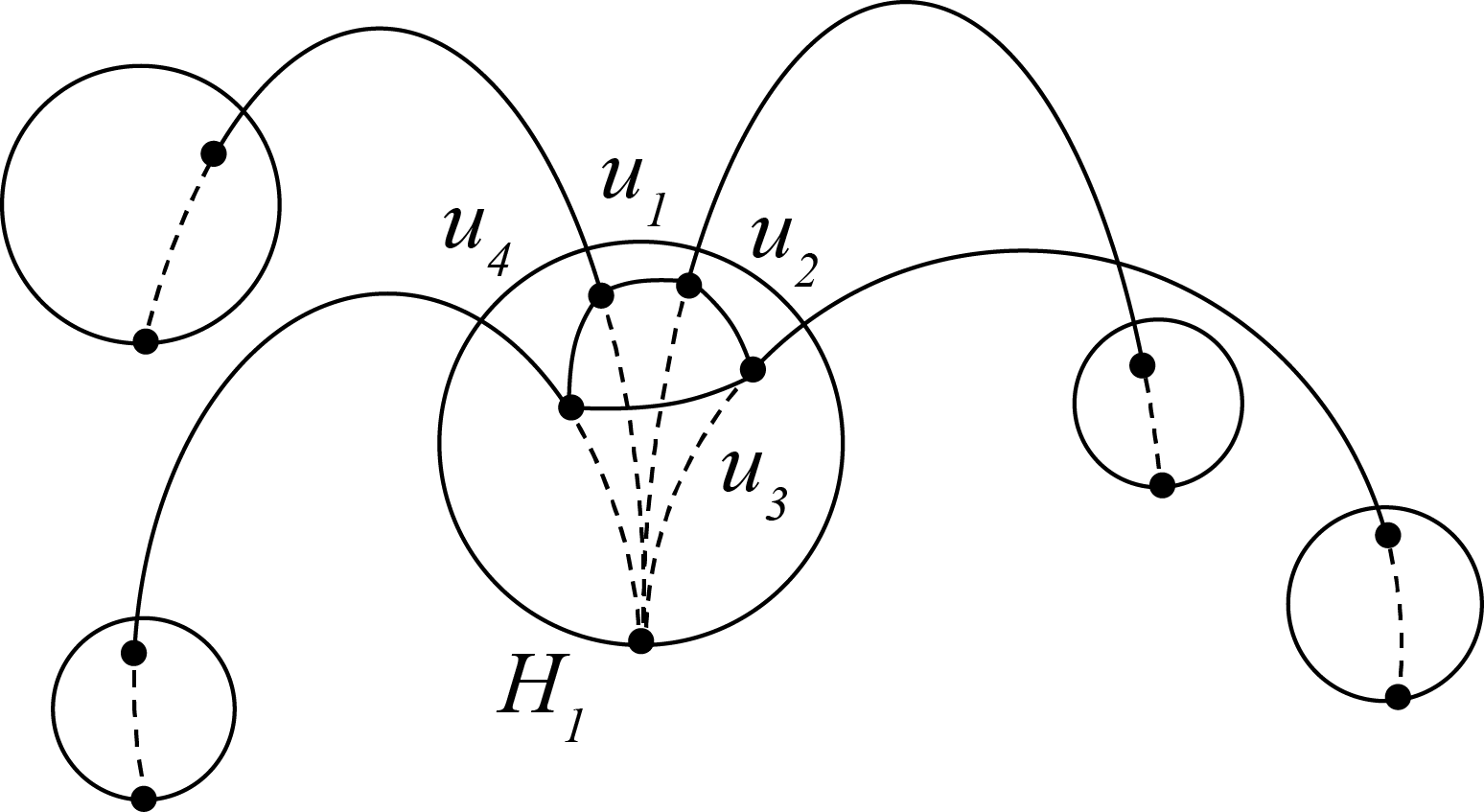} \hspace{0.3 in}
\includegraphics[scale=0.45]{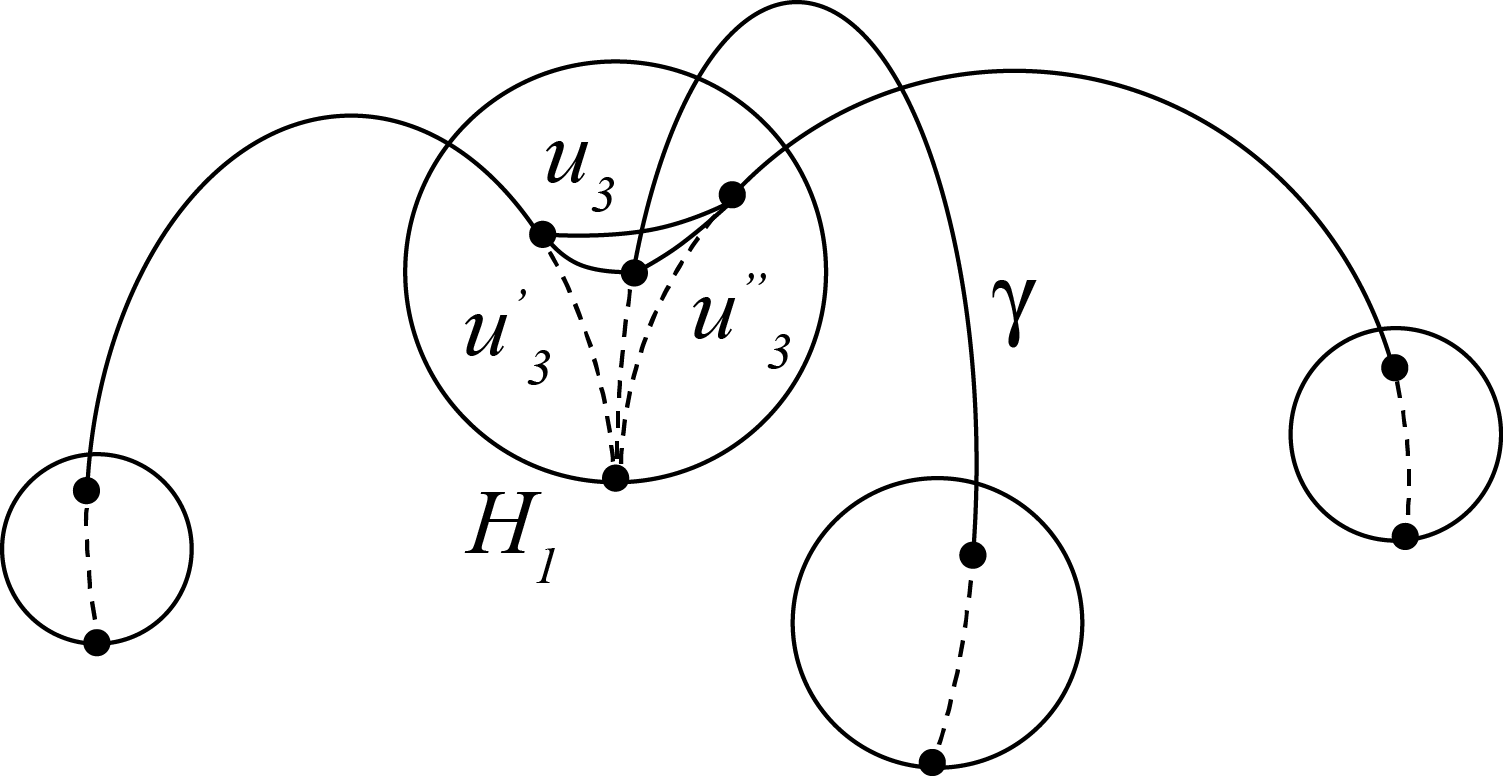}

Fig.11

\end{center}

Since the conditions (1)-(4) are satisfied by $\tau$, we may use Theorem 1.1 from \cite{PetronioWeeks} to state that the complete hyperbolic structure on $S^3-L$ exists. $\Box$

\vv
The diagram labels that satisfy the hyperbolicity equations and conditions (a)-(c) above were used to prove the existence of the complete hyperbolic structure on $L$. If, rather, we assume that the structure exists, we automatically obtain the following statement from our construction of a triangulation.

\vv

\textbf{Corollary 4.5.} Once a hyperbolic 3-manifold $M$ has a decomposition into two ideal cross-sectionally convex polyhedra in $\mathbb{H}^3$, there exists an ideal partially flat geodesic triangulation of $M$.

\section{Isotopy classes of crossing arcs}

The following statement provides a method for determining whether a crossing arcs is isotopic to a simple geodesic. Note that to check the conditions (a)-(c), one does not have to perform a decomposition of the link complement into polyhedra (or tetrahedra).
\vv

\textbf{Theorem 5.1.} Under the conditions (a)-(c) on the diagram labels of a reduced alternating diagram $D$ of a link $L$, its crossing arcs are isotopic to simple geodesics.

\vv

\textbf{Proof.} The proof of Theorem 1.1 from \cite{PetronioWeeks} implies that the topological space obtained by gluing the tetrahedra (some of which are possibly flat) is homeomorphic to $S^3-L$, and that the hyperbolic structure of $S^3-L$ locally induces its own metric on each tetrahedron. Therefore, the edges of $\tau$ are geodesics in $\mathbb{H}^3$, and the tetrahedra of $\tau$ are isometric to the corresponding tetrahedra in $S^3-L$, making the edges of the corresponding triangulation of $S^3-L$ geodesics as well. By construction, every crossing arc of $D$ is an edge of $\tau$.

 Let us now check that the crossing arcs are simple geodesics, \textit{i.e.} have no self-intersections. The only edges of $\tau$ that intersect in one point in $\mathbb{H}^3$ are two distinct edges of a flat tetrahedron (and, respectively, the edges that are identified with these two edges under the gluing). We need to check that these edges do not correspond to the same crossing arc in $S^3-L$. By the Thurston-Menasco construction, a crossing arc $c$ in $S^3-L$ corresponds to one edge $e_1$ of $\Pi_1$, and one edge $e_2$ of $\Pi_2$. Under the gluing, $e_1$ and $e_2$ are identified, and therefore they cannot be two distinct edges of one geometric tetrahedron intersecting in just one point. $\Box$ 

\section{Examples and conclusions}

The following examples illustrate that the conditions (a)-(c) are not too restrictive, and also show how to check them in practice. The first example is an infinite family of alternating braids, and the second one is an alternating link randomly chosen from the tables. Empirical data suggests that many other hyperbolic alternating links satisfy the conditions (a)-(c). 

\vv

\textbf{Example 6.1.} Consider an infinite family of alternating closed braids $(\sigma_1\sigma_3\sigma_2^{-1})^n$, $n>2$. Fig.12, left, shows a fragment of such a link with the diagram labels. The link can always be oriented as on the figure. Then the symmetry together with the relations for a 3-sided region imply that there are only two edge labels, $u_1$ and $u_2$, and two crossing labels, $w_1$ and $w_2$.

\begin{center}
\includegraphics[scale=1]{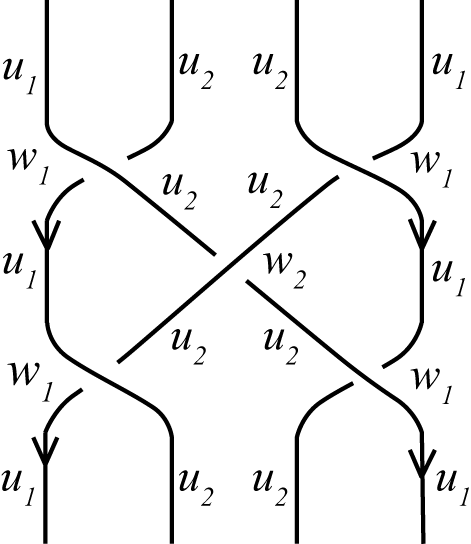} \hspace{0.5in}
\includegraphics[scale=0.6]{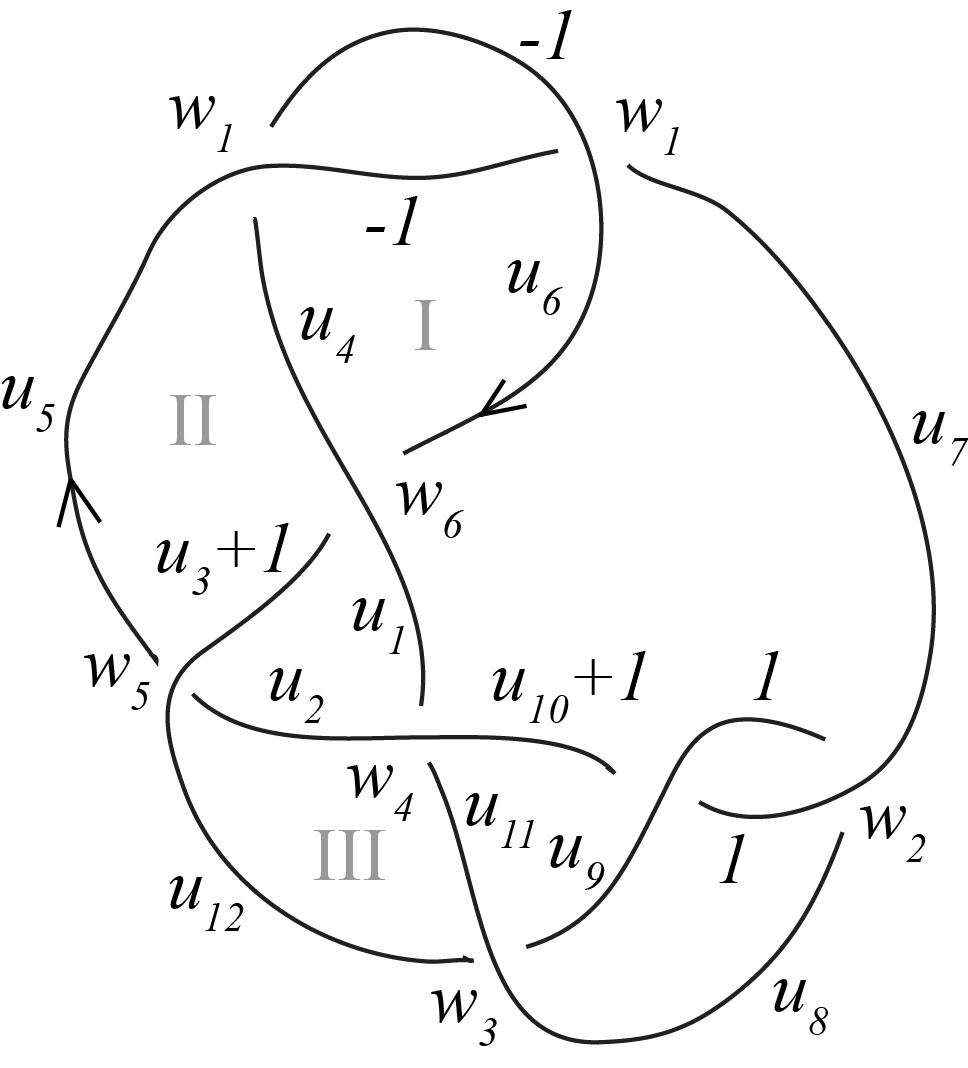}

Fig.12
\end{center}

The shape (in the sense of \cite{method}) of the regular $n$-sided region was established in \cite{method} and is \h $\frac{w_1}{u_1^2}=L=\frac{1}{4}\sec^2\frac{\pi}{n}$. \h Additionally, the labels satisfy the equations $-w_1=(u_1+1)(u_2+1)$, \h $\frac{w_1}{u_2^2}-\frac{w_2}{u_2^2}=1$,\h $w_2=-(u_2+1)^2$. Then $u_1$ satisfies $Lu_1^3-Lu_1^2-1=0$, and $u_2=\frac{Lu_1^2+1}{-2}$. 

A quick computation (one can use a computer algebra system) shows that $\operatorname{Im}u_1<0$. In particular, as $n$ approaches $\infty$, $\operatorname{Im}u_1$ monotonically decreases towards the limit of \h $\frac{2^{4/3} - (31 + 3 \sqrt{105})^{2/3}}{2^{
 2/3} \sqrt{3} (31 + 3 \sqrt{105})^{1/3}}$, which is approximately $-1.20563$, starting with a negative number for $n=2$. One can also check that $\operatorname{Im}u_2<0$. In fact, $\operatorname{Im}u_2$ monotonically increases towards the limit of

 \begin{center}
 $\frac{44 \sqrt{3} + 
  12 \sqrt{35} - (13 \sqrt{3} + 3 \sqrt{35}) (62 + 6 \sqrt{105})^{
   1/3}}{24 (62 + 6 \sqrt{105})^{2/3}}$, 
   \end{center}

    \noindent which is approximately $-0.0895648$, as $n \rightarrow \infty$. And $\operatorname{Im}{\frac{-u_1}{u_2+1}>0}$, $\operatorname{Im}{\frac{-u_2}{u_1+1}>0}$. Similarly, these are monotonic functions of $n$, and one can write down the positive numbers they converge to as $n\rightarrow \infty$. Lastly, $\frac{-w_1}{(u_1+1)(u_2+1)}$ is not purely real for all $n=4, 5, ...$. Therefore, the conditions (a)-(c) are satisfied.

A similar argument can be applied to any alternating closed braid of the type
\begin{center}
$(\sigma_1\sigma_3...\sigma_{2k+1}\sigma_2\sigma_4...\sigma_{2k}^{-1})^n$, $n>2$,
\end{center}  showing that the crossing arcs of the reduced alternating diagram of such a link are isotopic to geodesics. Note that these are braids that have an even index. One may also consider closed alternating braids with no bigons of odd index, generalizing Example 6.2 from \cite{method}.

\vv

\ind \textbf{Example 6.2.} Consider the 2-component link $8_8^2$ in the Rolfsen table, and its reduced alternating diagram (Fig.12, right). Assign diagram labels and orientation to it as on the figure. Recall that every region yields three region relations. Below we give the relations, and the decimal values of the labels necessary to check the conditions (a)-(c). This simple calculation shows that all the crossing arcs are isotopic to geodesics.

\v

 From the regions I, II, III respectively: 
 
\noindent  $w_1+u_6=0$, \h\h $w_6-u_6u_4=0$, \h\h $w_1+u_4=0$, \h\h $w_1+(u_4+1)(u_5+1)=0$, \h\h

\noindent $w_5-(u_5+1)(u_3+1)=0$, \h\h $w_6+(u_4+1)(u_3+1)=0$, \h\h $w_5-(u_{12}+1)(u_2+1)=0$, \h\h 

\noindent $w_4-(u_2+1)(u_{11}+1)=0$, \h\h $w_3-(u_{11}+1)(u_{12}+1)=0$. 

From the other three 3-sided regions: 

\noindent $w_6-u_3u_1=0$, \h\h $w_4+u_2u_1=0$, \h\h $w_5+u_2u_3=0$, \h\h $w_2-u_{10}u_9=0$, \h\h $w_3+u_{11}u_9=0$, \h\h $w_4+u_{10}u_{11}=0$, \h\h $w_2+u_8+1=0, \h\h w_3-(u_8+1)(u_9+1)=0, \h\h w_2+u_9+1=0$.

\v
 
 Lastly, of the two 5-sided regions (one of which is the outer region of the diagram), each yields equations of the form \h\h $\xi_1\xi_3-(\xi_1+\xi_2+\xi_3)+1=0$, \h\h $\xi_2\xi_4-(\xi_2+\xi_3+\xi_4)+1=0$, \h\h $\xi_3\xi_5-(\xi_3+\xi_4+\xi_5)+1=0$, \h\h where for the inner region \h\h $\xi_1=\frac{-w_1}{(u_6+1)(u_7+1)}$, \h\h  $\xi_2=\frac{-w_2}{u_7+1},$ \h\h $\xi_3=\frac{-w_2}{u_{10}+1},$ \h\h $\xi_4=\frac{w_4}{(u_1+1)(u_{10}+1)},$ \h\h $\xi_5=\frac{-w_6}{(u_1+1)(u_6+1)}$, \h\h and for the outer \h\h $\xi_1=\frac{-w_1}{u_7},$ \h\h $\xi_2=\frac{w_2}{u_7u_8},$ \h\h $\xi_3=\frac{-w_3}{u_{12}u_8},$ \h\h $\xi_4=\frac{-w_5}{u_{12}u_5},$ \h\h $\xi_5=\frac{w_1}{-u_5}$.
 
 \v
 
 One can then use a computer algebra system to obtain the solutions. The solution that satisfies the condition (a)-(c) has the following approximate decimal values of the labels: \h\h $w_1=0.37-0.52i$, \h\h $w_2=-0.37-0.52i$, \h\h $w_3=-0.13+0.39i$, \h\h $w_4=0.19+0.34i$, \h\h $w_5=-0.19+0.34i$,\h\h 
 
\noindent  $w_6=-0.13-0.39i,$ \h\h $u_1=-0.08+0.63i,$ \h\h $u_2=-0.5+0.36i,$ \h\h $u_3=-0.58+0.27i,$ \h\h 

\noindent $u_4=u_6=-0.37+0.52i,$ \h\h $u_5=-0.85+0.78i$, \h\h $u_7=-0.5+1.9i,$ \h\h $u_8=u_9=-0.63+0.52i$, \h\h $u_{10}=-0.05+0.78i,$ \h\h $u_{11}=-0.42+0.27i$, \h\h $u_{12}=-0.92+0.63i$. \h\h Therefore, the crossing arcs of the diagram on Fig.12, right, are isotopic to geodesics.
 \vv
 
It is also of note that experimental data suggests the following conjecture.

 \vv

\textbf{Conjecture 6.3.} The Thurston-Menasco ideal polyhedra in a hyperbolic alternating link complement in $S^3$ are cross-sectionally convex hyperbolic polyhedra, and the Thurston-Menasco ideal triangulation is a partially flat geodesic triangulation. 

\section{Acknowledgments }

I am grateful to Marc Lackenby for many enlightening conversations and comments on the draft; to Marc Culler, David Futer, Joel Hass, Carlo Petronio, Jessica Purcell, Morwen Thistlethwaite, and Abigail Thompson for helpful discussions. The work was partially supported by the NSF-AWM travel mentoring grant, and by the NSF DMS-1406588 grant.  

\vvv

\vv

\parbox[t]{4in}{
Anastasiia Tsvietkova\\
Department of Mathematics and Computer Science\\
Rutgers University, Newark\\
360 Dr. Martin Luther King Jr. Blvd.\\
Hill Hall 325, Newark, NJ 07102\\
a.tsviet@rutgers.edu
}


\begin{thebibliography}{90}



\bibitem{Adams} Colin C. Adams, \textit{Unknotting tunnels in hyperbolic 3-manifolds}, Math. Ann. 302 (1995), no. 1, 177--195.


\bibitem{AdamsReid} Colin C. Adams, Alan W. Reid, \textit{Unknotting tunnels in two-bridge knot and link complements}, Comment. Math. Helv. 71 (1996), no. 4, 617--627.


\bibitem{Akiyoshi}  H. Akiyoshi, M. Sakuma, M. Wada, Y. Yamashita, \textit{Ford domains of punctured torus groups and two-bridge knot groups}, Hyperbolic spaces and related topics, II (Japanese) (Kyoto, 1999). Sūrikaisekikenkyūsho Kōkyūroku No. 1163 (2000), 67--77.

\bibitem{Burton} S. D. Burton, J. S. Purcell, \textit{Geodesic systems of tunnels in hyperbolic 3-manifolds},  Algebr. Geom. Topol. 14 (2014), no. 2, 925--952.

\bibitem{Cooper} D. Cooper, D. Futer, J. Purcell, \textit{Dehn filling and the geometry of unknotting tunnels},
Geom. Topol., Vol. 17 (2013), no. 3, 1815-1876.



\bibitem{FuterAppendix} F. Gueritaud, \textit{On canonical triangulations of once-punctured torus bundles and two-bridge link complements}, with an appendix by D. Futer, Geom. Topol. 10 (2006), 1239--1284. 

\bibitem{Gueritad} F. Gueritaud, \textit{Geometrie hyperbolique effective et triangulations ideales cononiques en dimension 3}, Ph.D. thesis, L'Universit\'e de Paris (2006)


\bibitem{Lackenby} M. Lackenby, \textit{The volume of hyperbolic alternating link complements}, with an appendix by I. Agol and D. Thurston., Proc. London Math. Soc. (3) 88 (2004), no. 1, 204--224.


\bibitem{LackenbyHeegaard} M. Lackenby, \textit{An algorithm to determine the Heegaard genus
of simple 3–manifolds with nonempty boundary}, Algebr. Geom. Topol. 8 (2008) 911--934.



\bibitem{Menasco} W. W. Menasco, \textit{Polyhedra representation of link complements}, Low-dimensional topology (San Francisco, Calif., 1981), 305--325, Contemp. Math., 20, Amer. Math. Soc., Providence, RI, 1983.


\bibitem{NeumannZagier} W. Neumann, D. Zagier, \textit{Volumes of hyperbolic three-manifolds}, Topology 24 (1985), no. 3,
307--332.

\bibitem{Carlo} C. Petronio, \textit{An algorithm producing hyperbolicity equations for a link complement in $S^3$}, Geom. Dedicata 44 (1992), no. 1, 67--104.


\bibitem{PetronioWeeks} C. Petronio, J. R. Weeks, \textit{Partially flat ideal triangulations of cusped
hyperbolic 3-manifolds}, Osaka J. Math. 37 (2000), 453--466.


\bibitem{SakumaWeeks} M. Sakuma and J. R.
Weeks, \textit{Examples of canonical decomposition of hyperbolic link complements}, Japan. J. Math. (N. S.) 21 (1995), No. 2, 393--439.

\bibitem{Thurston} W. P. Thurston, \textit{The geometry and Topology of Three-Manifolds}, Electronic Version 1.1 (March 2002), http://www.msri.org/publications/books/gt3m/  

\bibitem{method} M. Thistlethwaite, A. Tsvietkova, \textit{An alternative approach to hyperbolic structures on link complements}, Algebr. Geom. Topol. 14 (2014), 1307--1337.

\bibitem{thesis} A. Tsvietkova, \textit{Hyperbolic Structures from Link Diagrams}, Ph.D. thesis, University of Tennessee (2012).


\bibitem{Weeks}
J. R. Weeks, \textsl{SnapPea}: a computer program for
creating and studying hyperbolic $3$--manifolds, freely available
from


http://thames.northnet.org/weeks/index/SnapPea.html



\end{thebibliography}
\end{document}